\documentclass{amsart}
\usepackage{amssymb}
\usepackage{graphicx}
\usepackage{epstopdf}
\newtheorem{thm}{Theorem}[section]
\newtheorem{cor}[thm]{Corollary}
\newtheorem{lem}[thm]{Lemma}
\newtheorem{lemdef}[thm]{Lemma and Definition}
\newtheorem{prop}[thm]{Proposition}
\newtheorem{conj}[thm]{Conjecture}
\theoremstyle{definition}
\newtheorem{defn}[thm]{Definition}
\theoremstyle{remark}
\newtheorem{rem}[thm]{Remark}
\numberwithin{equation}{section}
\newcommand{\To}{\longrightarrow}

\newcommand{\inv}{^{-1}}
\newcommand{\C}{\mathbb C}
\newcommand{\Z}{\mathbb Z}
\newcommand{\R}{\mathbb R}

\newcommand{\x}{\times}

\newcommand{\si}{\sigma}

\begin{document}

\title[A mirror construction and total positivity]
{A mirror construction for the 
totally nonnegative part of the
Peterson variety}
\author{Konstanze Rietsch}%
\address{King's College London, UK}%
\email{konstanze.rietsch@kcl.ac.uk}%


\thanks{
The author was supported by a Royal Society Dorothy Hodgkin
Research Fellowship. This work was additionally supported 
by FWF Project P~17108.}


\dedicatory{Dedicated to Professor George Lusztig on his 60th birthday}
\begin{abstract}
We explain how A. Givental's mirror symmetric family
\cite{Giv:MSFlag} to the type $A$ flag variety and its
proposed 
generalization \cite{BCKS:MSPFlag} to partial flag varieties
by Batyrev, Ciocan-Fontanine, Kim and van Straten 
relate to the Peterson variety $Y\subset SL_{n}/B$. 
We then use this theory to describe the totally
nonnegative part of $Y$, extending a result from \cite{Rie:QCohPFl}.
\end{abstract}
\maketitle
\section{Introduction}
The (type $A$) Peterson variety is a
remarkable $(n-1)$-dimensional 
projective subvariety
of the full flag variety $SL_{n}/B$ used by 
Dale Peterson to construct all of the 
small quantum cohomology rings of 
the partial flag varieties
$SL_{n}/P$. This paper has two aims:
firstly 
to relate the mirror symmetry constructions
of Givental \cite{Giv:MSFlag} 
and Batyrev, Ciocan-Fontanine, Kim, and van Straten
\cite{BCKS:MSPFlag} to the Peterson variety,
and secondly to use these constructions 
to describe the totally nonnegative
part of the Peterson variety.

The mirror constructions of  \cite{Giv:MSFlag} and \cite{BCKS:MSPFlag}  
provide in the full flag variety case, and provide conjecturally in the partial flag variety case, 
a set of solutions to the quantum cohomology 
$D$-module 
-- a system of 
differential equations introduced by Givental whose 
total symbols recover relations of the 
small quantum cohomology ring \cite{Giv:EquivGW} --
in terms of oscillating integrals along families of cycles lying in a 
`mirror family'. These mirror families
are $k$-dimensional families $Z\to \C^k$ of
affine varieties of the same dimension as $SL_n/P$
which are defined in terms of an associated graph, 
see Figure~\ref{f:graph}, and which come with natural 
volume forms on the fibers and a
phase function $\mathcal F:Z\to \C$. 
Here $k=\dim H^2(SL_n/P)$.
If the mirror conjecture holds then critical points of the 
phase function $\mathcal F$
along individual mirrors should relate to elements in  
the spectrum of the quantum cohomology ring, where fixing a variety in the mirror family 
corresponds to fixing the values of the quantum parameters $q_1,
\dotsc, q_k$ in  $qH^*(SL_n/P)$.

In Dale Peterson's theory the spectrum of 
$qH^*(SL_{n}/P)$ is precisely
a stratum $Y_P$ of the Peterson variety $Y$. 
In this paper we compare Peterson's $Y_P$ with the variety $Z^{crit}$ 
swept out by the critical points of $\mathcal F$ along the fibers of the mirror 
family $Z$ from \cite{BCKS:MSPFlag, Giv:MSFlag}. 
As it turns out 
$Z^{crit}$ recovers the parts of the Peterson variety that 
lie in certain Deodhar strata (a finer decomposition of the flag variety
than the Bruhat decomposition). 
If $P=B$ then this includes an open dense subset 
of $Y_B$. But in the case $P\ne B$, the variety
$Y_P$ can have entire irreducible components which 
lie in the `wrong' Deodhar stratum and hence are not seen by $Z$. 
This phenomenon is demonstrated explicitly in 
Section~\ref{s:Gr2(4)} for $SL_4/P=Gr_2(\C^4)$. 

In this special case, $Gr_2(\C^4)$,
an earlier mirror construction consistent with the `GBCKS' mirror construction from 
\cite{Giv:MSFlag,BCKS:MSPFlag} was 
given by Eguchi, Hori and 
Xiong in  \cite[Appendix B]{EHX:GravQCoh}, see \cite{BCKS:MSGrass}. 
Its deficiency with regard to recovering the 
quantum cohomology ring was observed also in \cite{EHX:GravQCoh}, where it was
fixed in an ad hoc way by a partial compactification. 
For a `fix' of the GBCKS construction for general $SL_n/P$
we refer to our sequel paper \cite{Rie:MSgen}. It has not been checked how in the case 
of $Gr_2(\mathbb C^4)$ the general construction of \cite{Rie:MSgen} 
relates to the ad hoc construction from \cite{EHX:GravQCoh}.

Next we turn our attention to total positivity. The totally nonnegative part 
$(SL_{n}/B)_{\ge 0}$ of the 
flag variety was defined by Lusztig~\cite{Lus:TotPos94}
as an extension of the classical theory of total
positivity for matrices. It is a semi-algebraic subset
inside the real flag variety $SL_n(\R)/B$ (which we view 
with its Hausdorff topology).

In \cite{Rie:QCohPFl} we showed that the totally positive
part of $Y_P$ (that is, the open interior of $Y_P\cap (SL_n/B)_{\ge 0}$)
agrees with the subset of $Y_P$ where all of the Schubert
classes take positive real values. Using this result it was then proved
that the quantum parameters restrict to give a homeomorphism
$
Y_{P,>0}\overset\sim\To \R_{>0}^k,
$
where $k=\dim Y_P$, making $Y_{P,>0}$ a cell.

In Section~\ref{s:TotPos} we use the mirror constructions  
from the previous sections to give a direct new proof of the
above parameterization. In fact we can extend the result to the boundary
to get a homeomorphism,
$$
Y_{P,\ge 0}\overset\sim\To \R_{\ge 0}^{k},
$$ 
parameterizing the totally nonnegative part of $Y_P$. 
Therefore we obtain a cell decomposition 
of the whole totally nonnegative part of the Peterson variety 
$Y$. This mirror symmetric approach to proving the cell decomposition 
has the advantage of being
completely elementary, whereas 
the proof in \cite{Rie:QCohPFl} relied 
on positivity of the  
structure constants 
of the quantum cohomology  rings involved (the 
3-point genus zero Gromov-Witten invariants for 
$SL_{n}/P$). On the other hand, though, 
we obtain no results
about positivity of Schubert classes
using only the mirror construction.  

Finally, it is shown that the totally nonnegative
part $Y_{\ge 0}$ of the Peterson variety is contractible. 
We conjecture that $Y_{\ge 0}$, as a
cell decomposed space, is homeomorphic to an
$(n-1)$-dimensional cube. 

The interpretation of the GBCKS mirror 
construction and resulting proof of the
cell decomposition of the totally nonnegative part of the Peterson variety $Y$ 
presented here date back to 2002, and 
were presented at the Erwin Schroedinger Institute in 
January of 2003 as well as alluded to
in a footnote in \cite{Rie:EWM}. In the full flag variety case a similar interpretation 
(but very different application) of Givental's mirror coordinates has since appeared 
also in the interesting work of 
Gerasimov, Kharchev, Lebedev and Oblezin \cite{GKLO:GaussGivental}
on the quantum Toda lattice.   
\vskip 3mm
\noindent{\sl Acknowledgements~:} I would particularly like to thank 
George Lusztig
and Dale Peterson. The first for introducing me to the marvelous theory
of total positivity,   
and the second for his inspiring lectures on quantum 
cohomology. Without either one of them this paper would not 
have been written. These results were mostly written up while on leave in
 Waterloo, Canada.  I thank the University of Waterloo for its hospitality. 

\section{Notation}

From now on we let $n$ be the rank. 
Consider $G=SL_{n+1}(\C)$
with fixed Borel subgroups $B=B^+$, the group of
upper-triangular matrices, and $B^-$ the lower-triangular
matrices, and with maximal torus $T=B^+\cap B^-$. We also 
have $U^+$ and $U^-$, the unipotent radicals of $B^+$ and $B^-$,
respectively.
Let $I=\{1,\dotsc, n\}$ and $e_i,f_i$ 
the usual Chevalley generators of the Lie algebra 
$\mathfrak g=\mathfrak{sl}_{n+1}$. So $e_i$ is the matrix
with $1$ in position $(i, i+1)$ and $0$ everywhere else, and $f_i$ 
is its transpose. Let  
\begin{equation*}
x_i(t):=\exp(t e_i),\quad y_i(t):=\exp(t f_i), \quad t\in\C
\end{equation*}
be the associated simple root subgroups. The datum 
$(T,B^+,B^-,x_i,y_i;i\in I)$ is called a pinning by Lusztig
\cite{Lus:TotPos94}.

The Weyl group 
$W=N_G(T)/T$ is isomorphic to the symmetric group 
$S_{n+1}$. Define
representatives
\begin{equation*}
\dot s_i:=y_i(-1)x_i(1)y_i(-1),\quad i\in I.
\end{equation*} 
for the simple reflections $s_i:=\dot s_i T$. The $s_i$  are 
Coxeter generators for $W$. For general $w\in W$ a representative
$\dot w\in G$ can be defined by $\dot w=\dot s_{i_1}\dot
s_{i_2}\cdots \dot s_{i_m}$, where $ s_{i_1} s_{i_2}\cdots
s_{i_m}$ is a (any) reduced expression for $w$. The length $m$ of
a reduced expression for $w$ is denoted by $\ell(w)$.

Let $P\supseteq B$ be a  parabolic subgroup of $G$. Then
there is a corresponding parabolic subgroup $W_P$ of $W$
generated by the elements $s_i$ with $\dot s_i\in P$. Define 
$I_P=\{i\in I\ |\ \dot s_i\in P \}$ and $I^P$ its complement in
$I$. We have
\begin{align*}
 W_P&=\left<s_i\ |\ i\in I_P\right >, \\
 W^P&:=\{w\in W\ |\ \ell(w s_i)>\ell(w) \text{ for all $i\in I_P$}.
 \}
\end{align*}
The longest element in $W_P$ is denoted by $w_P$. The longest 
element in $W$ is also denoted $w_0$.  

Let $I^P=\{n_1,\dotsc,n_k\}$ where
$0=n_0<n_1<n_2<\cdots<n_k<n+1=n_{k+1}$. Then the homogeneous space
$G/P$ can be identified with the variety of partial flags
 $$
\mathcal F_{n_1,n_2,\dotsc, n_k}(\C^{n+1})= \{\,\{0\}\subset
V_1\subset V_2\subset\dots\subset V_{k}\subset \C^{n+1}\, |\,
\dim_\C(V_j)=n_j\, \}.
 $$

\section{Quantum cohomology of $SL_{n+1}/P$}\label{s:qcoh}

Let $H^*(G/P):=\bigoplus_k H^{2k}(G/P)$ be the cohomology of
$G/P$ viewed as a graded vector space with grading given by $k$.
We will always take coefficients in $\C$. 
For $w\in W^P$ denote by $\sigma^{w}_P$ the Poincar\'e dual
class to the Schubert cycle $[X_w]$ where $X_w=\overline{B^- w P/P}$. 
It is well known that the Schubert classes $\sigma^{w}_P$
are a homogeneous basis of $H^*(G/P)$ with $\deg(\sigma_P^w)=\ell(w)$.

The small quantum cohomology ring of the partial flag variety
$SL_{n+1}/P$ has been described in the papers
\cite{AstSa:QCohPFl,Cio:QCohPFl,Kim:QCohPFl}. As a graded vector
space it is given by
 $$
 qH^*(SL_{n+1}/P)=H^*(SL_{n+1}/P)\otimes\C[q_1^{P},\dotsc, q_k^{P}],
 $$
where $\C[q_1^{P},\dotsc, q_k^{P}]$ is a graded polynomial ring
with $\deg (q^P_j)=n_{j+1}-n_{j-1}$. The multiplicative structure
constants are $3$-point genus~$0$ Gromov-Witten invariants, see
for example \cite{GiKi:FlTod,Cio:QCohPFl,FoGePo:QSchub,Kos:QCoh}
or \cite{CoxKatz:QCohBook,McDSal:QCohBook}.
For the purposes of this paper we will
be  mainly interested in 
presentations of these rings. 

\subsection{}\label{s:presentations}
Let  
$$\C[\mathfrak h]=
Sym^\bullet(\mathfrak h^*)=\C[x_1,\dotsc, x_{n+1}]/(x_1+\dotsc + x_{n+1})
$$ 
be the coordinate ring of $\mathfrak h=Lie(T)$, 
where the $x_i$ are the coordinates corresponding to 
the matrix entries along the diagonal. 
The $\Z$-span of
the $x_i$ is the character lattice $X^*(T)$ inside $\mathfrak h^*$. 
The assignment taking  a character $\lambda $ to the first Chern
class of the associated line bundle $\mathcal L_\lambda =G\x_B
\C_\lambda$ on $G/B$, extends to a ring homomorphism
$\C[\mathfrak h]\to H^*(G/B)$. 
By Borel \cite{Borel:CohG/P}, this map identifies
$H^*(G/B)$ with the quotient  
$$
\C[x_1,\dotsc,
x_{n+1}]/(e_1^{(n+1)},\dotsc, e_{n+1}^{(n+1)}),
$$ where
$e_l^{(n+1)}=e_l(x_1,\dotsc, x_{n+1})$ is the $l$-th elementary
symmetric polynomial in $n+1$ variables.
Moreover the projection $G/B\to G/P$ gives rise to an inclusion  
$H^*(G/P)\to H^*(G/B)$ which identifies $H^*(G/P)$ with the
$W_P$-invariant part of $H^*(G/B)$. 
Explicitly, consider the ring $\C[x_1,\dotsc,x_{n+1}]^{W_P}$,
which is a polynomial ring generated by the elementary symmetric polynomials 
\begin{eqnarray*}
&\sigma^{(1)}_l:=e_l(x_1,\dotsc,x_{n_1}), &\qquad l=1,\dotsc,n_1,\\
&\sigma^{(2)}_l:=e_l(x_{n_1+1},\dotsc, x_{n_2}), 
&\qquad l=1,\dotsc , n_2-n_1,\\
&\qquad\qquad\vdots&\\
&\sigma^{(k+1)}_l:=e_l(x_{n_k+1},\dotsc, x_{n+1}),& \qquad l=1,\dotsc,
n+1-n_k.
\end{eqnarray*}
The full elementary symmetric polynomials $e_r^{(n+1)}$ may 
be expressed as polynomials in the $\sigma^{(j)}_l$ and we let $J$ denote
the ideal these polynomials generate. 
Then we have
\begin{equation}\label{e:BorelIso}
H^*(G/P)\cong\C[\sigma^{(1)}_1,\sigma^{(1)}_2,\dotsc,\sigma^{(k+1)}_{n+1-n_k}]/J.
\end{equation}

\subsection{}\label{s:qpresentations}
The analogous presentation 
of the quantum cohomology ring due 
to \cite{AstSa:QCohPFl,Cio:QCohPFl,Kim:QCohPFl} goes as follows.
From now on let us write $\sigma^{(j)}_l$ for the element
$\sigma^{(j)}_l\otimes 1\in qH^*(G/P)$, and similarly $q_j^P$ or
just $q_j$ for $1\otimes q_j^P$. These are the generators.
\begin{defn}[$(\mathbf q,P)$-elementary symmetric polynomials]
\label{d:E} Let $l\in\Z$ and $j\in\{-1,0,\dotsc,k+1\}$. Define
elements
$E^{(j)}_{l,P}=E^{(j)}_{l}\in\C[\si^{(1)}_1,\dotsc,\si^{(k+1)}_{n+1-n_k}\,
,\,  q_1,\dotsc,q_k]$ recursively as follows. The initial values
are
\begin{equation*}
E^{(-1)}_l= E^{(0)}_l=0\ \text{ for all $l$},\  \ \text{and }\
E^{(j)}_{l}=0 \ \text{ unless } 0\le l\le n_j ,
\end{equation*}
and we set $\si^{(j)}_l=0$ if $l>n_j-n_{j-1}$ and
$\sigma^{(j)}_0=1$ for all $j$. For $1\le j\le k+1$ and $0\le l\le n_l$ the
polynomial $E^{(j)}_l$ satisfies
\begin{multline*}
E^{(j)}_l\hskip -.1cm=\hskip -.1cm\left (E^{(j-1)}_{l} +
\si^{(j)}_1 E^{(j-1)}_{l-1}+\cdots+ \si^{(j)}_{l-1} E^{(j-1)}_{1}+
\si^{(j)}_{l} \right ) \\ +(-1)^{n_{j}-n_{j-1}+1}q_{j-1}
E^{(j-2)}_{l-n_{j}+n_{j-2}}.
\end{multline*}
\end{defn}

\begin{thm}[\cite{AstSa:QCohPFl,Kim:QCohPFl,Cio:QCohPFl}]
\label{t:presentationP} The quantum cohomology ring $ qH^*(G/P)$
is given by the generators
$\sigma^{(1)}_1,\dotsc,\sigma^{(k+1)}_{n+1-n_k}, q_1,\dots, q_{k}$
with relations
\begin{equation*}
 E^{(k+1)}_1= E^{(k+1)}_2=\cdots = E^{(k+1)}_{n+1}=0.
\end{equation*}
\end{thm}

\section{The Peterson variety}

Dale Peterson \cite{Pet:QCoh} 
discovered a remarkable unified construction for all of 
the quantum cohomology rings $qH^*(G/P)$, for varying $P$,
as coordinate rings of the strata of a single projective variety $Y$.
For $G$ of general type this 
`Peterson variety' $Y$ is a subvariety of the  
Langlands dual flag variety $G^\vee/B^\vee$. 
We will recall his result in type $A$.

\subsection{}
In our conventions the Peterson variety will be 
a subvariety of $G/B^-$, where $G=SL_{n+1}(\C)$. Let us recall 
first the
Bruhat and opposite Bruhat decompositions
\begin{equation*}
G/B^-=\bigsqcup_{w\in W} B^-\dot w B^-/B^- =
\bigsqcup_{v\in W} B^+\dot v B^-/B^-.
\end{equation*}
We also define
\begin{equation*}
\mathcal R_{v,w}:=B^+\dot v B^-\cap B^-\dot w B^-/B^-.
\end{equation*} 
This intersection of opposed Bruhat cells is 
smooth of pure dimension $\ell(w)-\ell(v)$ if 
$v\le w$ in the Bruhat order, and otherwise empty,
see \cite{KaLus:Hecke2,Lus:IntroTotPos}.

Let $\{\omega_i\ |\ i\in I\}$ be the set of fundamental
weights. Consider $V^{\omega_r}=\bigwedge^r\C^{n+1}$,
 the $r$-th fundamental representation of $G$ with its standard basis 
$\{v_{i_1}\wedge\cdots\wedge v_{i_r}\,|\, 
1\le i_1<i_2\cdots <i_r\le n+1\}$. The stabilizer of the 
highest weight space $\left< v_1\wedge\cdots\wedge v_r\right>_\C$
defines a maximal parabolic which we denote $P_{\omega_r}$.
Let us write
$V_{-\omega_r}$ for
$V^{\omega_{n-r+1}}$, which is the representation with lowest weight $-\omega_r$,
and fix the lowest
weight vector $v_{-\omega_{r}}=v_{r+1}\wedge\cdots\wedge v_{n+1} $.
For
$w\in W^{P_{\omega_r}}$ we have a well defined rational function 
 \begin{equation}\label{e:M}
M_{w\omega_r}(g B^-):=
\frac{\left < g\cdot v_{-\omega_r},\dot w\cdot v_{-\omega_r}\right>}
{\left < g\cdot v_{-\omega_r},v_{-\omega_r}\right>}
 \end{equation}
on the flag variety $G/B^-$, where $\left<\ ,\ \right>$ 
denotes the inner product on $V_{-\omega_{r}}$ such that the standard
basis is orthonormal. 

Let us introduce the principal nilpotent element 
 \begin{equation*}
 \label{e:f}
f=f_1+\dotsc + f_{n}.
\end{equation*}
We write
$g\cdot X:=g X g\inv$ for the adjoint action of $g\in G$ on  
$X\in \mathfrak g$.

\begin{defn}[The Peterson variety]
Let $Y\subset G/B^-$ be the projective variety defined by
 $$
Y:=\left\{gB^- \,\left |\, g\inv\cdot f\in\mathfrak b^-\oplus
\sum_{i\in I}\C e_i\right.\right\}.
 $$  
More formally, $Y$ is defined by the equations 
$$
{pr}_{\mathfrak g_{\alpha}}(g\inv\cdot f)= 0,
$$
where ${pr}_{\mathfrak g_\alpha}$ is the 
projection onto the weight space $\mathfrak g_{\alpha}$, and
$\alpha$ runs through the 
set of all roots which are positive but not simple.  
For a parabolic $P\supseteq B$ define the (non-reduced) intersection
\begin{equation*}
Y_P:=Y\x_{G/B^-} B^+\dot w_P B^-/B^-.
\end{equation*}
Suppose $P'\supseteq P$ is another parabolic. Then we set
$$
Y_{(P,P')}:=Y\x_{G/B^-} \mathcal R_{w_P,w_{P'}}.
$$
We also write $Y_{P}^\circ$ for $Y_{(P,G)}$. 
\end{defn}

The Peterson variety and some generalized versions of it are also 
of independent interest and have 
been studied in the papers \cite{BriCar:PetVar,Kos:QCoh,Kos:QCoh2,Tym:paving}.

\subsection{} We now 
state Peterson's result 
in type $A$, see also \cite{Kos:QCoh} and \cite{Rie:QCohPFl,Rie:ErratumJAMS}.

\begin{thm}[Peterson \cite{Pet:QCoh}]\label{t:Peterson}
\begin{enumerate}
\item The $\C$-valued points of $Y$ decompose into a union of strata,
 $$
Y(\C)=\bigsqcup_{P\supseteq B} Y_P(\C).
 $$ 
\item
Let $w^{[r]}_l=s_{r-l+1} s_{r-l+2}\cdots s_{r-1}s_{r}$, where $1\le l\le r$.
For each parabolic $P$ there is a unique isomorphism
 $$
\psi_P:\C[Y_P]\overset\sim\To qH^*(G/P),
 $$
such that $M_{w^{[n_j]}_l \omega_{n_j}}\mapsto E^{(j)}_{l}$
for $j=1,\dotsc,k$ and $1\le l\le n_j$.
\item 
$\psi_P$ induces an isomorphism
\begin{equation*}
\psi_P^\circ:\C[Y_P^\circ]\overset\sim\To qH^*(G/P)[q_1\inv,\dotsc, q_k\inv].
\end{equation*} 
\end{enumerate}
\end{thm}
Note 
that $M_{w^{[r]}_l\omega_r}$ is a 
regular function on the Bruhat cell $B^+ \dot w_P B^-/B^-$ 
if $r\in I^P$ and $1\le l\le r$.

\section{The GBKCS mirror construction for $SL_{n+1}/P$}
In \cite{Giv:MSFlag},
A. Givental introduced a mirror family to the full flag variety 
$SL_{n+1}/B$ 
and proved a kind of mirror theorem.  
His mirror construction was generalized by Batyrev, Kim, Ciocan-Fontanine and 
van Straten in  \cite{BCKS:MSPFlag}, who defined a similar family 
associated to partial flag varieties 
$SL_{n+1}/P$ and conjectured the analogous
mirror theorem. We recall their construction, which we will refer to as
the GBCKS construction, here.

\subsection{}
Let us fix the 
partial flag variety 
$$SL_{n+1}/P=\mathcal F_{n_1,\dotsc, n_k}(\mathbb C^{n+1}).$$
As before $n_{k+1}=n+1$ and $n_0=0$. 
Define an oriented graph $(\mathcal V,\mathcal A)=
(\mathcal V^P,\mathcal A^P)$ as follows. Let the  
vertex set $\mathcal V^P\subset \Z^2$ be defined by  
$\mathcal V^P=
\mathcal V^P_{\star}\sqcup \mathcal V^P_{\bullet}$ where
\begin{align*}
\mathcal V^P_\bullet&=
\{(m,r)\in \mathbb Z_{\ge 0}^2\ |\ n_1\le m\le n, \text{ and $1\le r\le n_{j}$
if  $m< n_{j+1}$, for $j=1,\dotsc, k$}\},\\
\mathcal V^P_{\star}&=\{\star_j=(n_{j}-1,n_{j-1}+1)\ |\ j=1,\dotsc, k+1\}.
\end{align*}
Consider $v=({v_1},{v_2})$ in $\mathcal V^P$. If $v':=(v_1, v_2-1)$ is in 
$\mathcal V^P$ then there is a 
horizontal arrow, denoted $d_v$ or $d_{v_1,v_2}$,
pointing from $v$ to $v'$. If   
$v''=(v_1-1, v_2)$ is in $\mathcal V^P$ then 
there is a vertical arrow $c_{v}$, or $c_{v_1,v_2}$, 
going from $v$ to $v''$. We define $\mathcal A^P$ 
to be the 
set of all such arrows.

\begin{figure}
\begin{center}
\leavevmode
 \[ \includegraphics{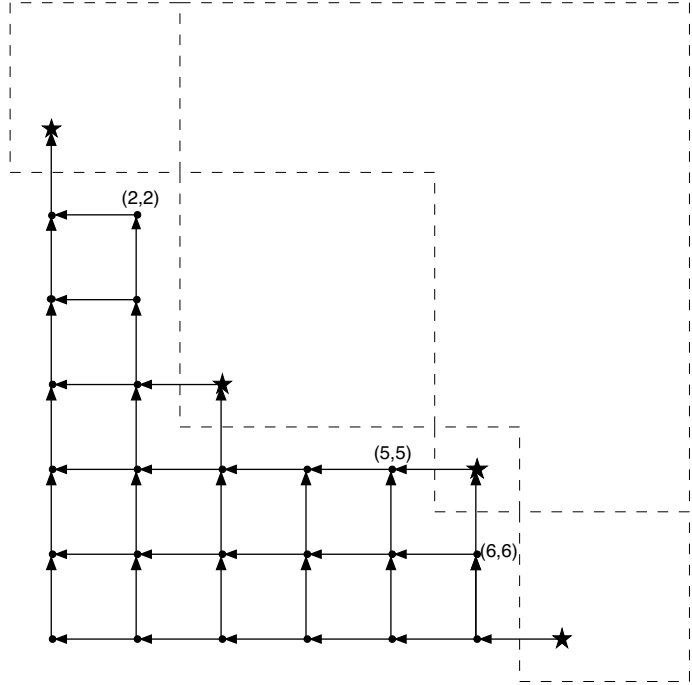}\]
\caption{\label{f:graph}The graph for $G/P=\mathcal F_{2,5,6}(\C^{8})$}
\end{center}
\end{figure}

See Figure~\ref{f:graph} for an example of a graph
$(\mathcal V^P,\mathcal A^P)$. The vertices are arranged like 
entries in a matrix, with a vertex $(i,j)$ positioned in the 
$i$-th row and $j$-th column. The 
dotted lines indicate the shape of the parabolic subgroup $P$
in question.
And the vertices in  $\mathcal V^P_{\star}$ and 
$\mathcal V^P_{\bullet}$ are represented 
by stars and dots, respectively. 

As the parabolic will be
fixed most of the time we may omit the superscript $P$ and write 
$\mathcal V$ for $\mathcal V^P$ and $\mathcal A$ for $\mathcal A^P$.

\subsection{}
Let
\begin{equation*}
Z=Z_P:=\left\{\rho=(\rho_a)_{a\in \mathcal A}\in \C^{\mathcal A}\ \left |\ 
\rho_c \rho_d=\rho_{d'}\rho_{c'},\ 
\begin{array}{l}
\text{whenever $c,c',d,d'$ form
}\\ \text {a square, \eqref{e:box}, in the graph}\  
\end{array}\right.
\right\}.
\end{equation*}
\begin{equation}\label{e:box}
\begin{matrix}\bullet &\overset{d}\longleftarrow&\hskip -.3cm \bullet\\
\hskip- .3cm c'\uparrow & &\uparrow c\\
\bullet &\overset{d'}\longleftarrow &\hskip -.3cm  \bullet
\end{matrix}
\end{equation}
The upper right hand corner vertex in \eqref{e:box} may 
of course lie in $\mathcal V_{\star}$. 

For simplicity of notation we identify the arrows with functions
on $Z$ via 
$$
a: \rho\mapsto \rho_a.
$$ 
The coordinate ring $\C[Z]$ 
can be viewed as the affine algebra over $\C$ 
with generators  $a\in\mathcal A$ and 
relations $c d=d'c'$ for $c,d,c',d'\in \mathcal A$ 
arranged as in \eqref{e:box}. We will refer to these
as `box relations'. There is a grading on 
$\C[Z]$ given by setting $deg(a)=1$ for 
every generator $a\in\mathcal  A$.

\subsection{}
The coordinate ring
$\C[Z]$ has some special elements which we define below. For 
$j=1,\dotsc k$ let $\tilde q_j$ be a product of generators represented by
the arrows along a path  
from vertex $\star_{j+1}$ to $\star_{j}$. Explicitly,  
\begin{equation*}
\tilde q_j= d_{n_{j+1}-1, n_{j+1}}
\left(\prod_{i=1}^{n_{j+1}-n_j-1} c_{n_j+i, n_j} \right)
\left(\prod_{i=1}^{n_j-n_{j-1}-1} d_{n_j,n_{j-1}+i+1}\right) 
c_{n_j, n_{j-1}+1},
\end{equation*}
where we have chosen the path along the outer rim. Note that 
$deg (\tilde q_j)=n_{j+1}-n_{j-1}$. 
Now $Z$ is viewed as a family of varieties via
\begin{equation}\label{e:family}
\tilde q=(\tilde q_1,\dots, \tilde q_{k}):Z\To \C^{k}.
\end{equation}
The fiber over $\tilde Q\in\C^{k}$ is denoted by $Z_{\tilde Q}$. 

\subsection{}
Finally \cite{Giv:MSFlag,BCKS:MSPFlag} introduce a function 
$$\mathcal F=\sum_{a\in \mathcal
A }a$$ 
on $Z$. This is the phase function of the proposed mirror model, see
Section~\ref{s:MS}. We will study
its critical point sets along the fibers of the family $Z$
in Section~\ref{s:main}.

\subsection{}\label{s:Zcirc}
Define
\begin{equation}
Z^\circ=Z^\circ_P:=
\left\{\rho\in Z\ |\ \rho_a\ne 0,\ \text{all $a\in\mathcal A$}\right\}=
\left\{\rho\in Z\ |\ \tilde q_j(\rho)\ne 0, \ 1\le j\le k\right\}.
\end{equation}
Let the map \eqref{e:family} restricted to $Z^\circ$ be again denoted by 
$\tilde q$,
\begin{equation*}
\tilde q=(\tilde q_1,\dots,\tilde  q_{k}):Z^\circ\To (\C^*)^{k}.
\end{equation*}
This restricted map is a trivial bundle with fiber isomorphic to
$(\C^*)^{\mathcal V_{\bullet}}$.

As in \cite{Giv:MSFlag} one can choose 
an explicit trivialization by introducing vertex 
variables $(t_v)_{v\in\mathcal V}$ running through $\C^*$.
For any arrow $a$
denote by $h(a)$ and $t(a)\in\mathcal V$ the head and tail of $a$.
Then $(t_v)_{v\in\mathcal V}
\mapsto \rho=(t_{h(a)}t_{t(a)}\inv)_{a\in\mathcal A}$
defines a map 
\begin{equation}\label{e:pre-trivialization}
(\C^*)^{\mathcal V}
\To Z^\circ.
\end{equation}
This map descends to the 
quotient by the diagonal action of $\C^*$
 to give an isomorphism
$(\C^*)^\mathcal V/ \C^*\overset\sim\To Z^\circ$.

Moreover, for given $\tilde Q=(\tilde Q_1,\dotsc,\tilde  Q_k)\in(\C^*)^k$, 
the map obtained from \eqref{e:pre-trivialization}
after fixing the $t_{\star_j}$ (uniquely up
to a common scalar multiple) such that 
$t_{\star_{j}}t_{\star_{j+1}}\inv=\tilde Q_j$ gives rise 
to an isomorphism
\begin{equation*}
(\C^*)^{\mathcal V_{\bullet}}\overset\sim\To Z_{\tilde Q}.
\end{equation*} 
Choosing $t_{\star_j}=\tilde Q_j\dotsc \tilde Q_{k}$, say, and $t_{\star_{k+1}}=1$
gives rise
to a global trivialization of $\tilde q:Z^\circ\to(\C^*)^k$. 
\subsection{}\label{s:Zparts}
For a pair of parabolics $P'\supseteq P$ containing $B$
the corresponding vertex sets are related by 
$\mathcal V^{P'}_\bullet\subseteq\mathcal V^{P}_{\bullet}$ and  we define
\begin{equation*}
Z_{(P,P')}:=\left\{\rho\in Z_P\ |\ \text{If $a\in\mathcal A^P$, then 
$a(\rho)= 0$ 
$\iff$
$h(a)\in \mathcal V^{P'}_\bullet$ or $t(a)\in\mathcal V^{P'}_\bullet$}
\right\}.
\end{equation*}
Note that if $P'=G$ we have $Z_{(P,G)}=Z^\circ_P$. In general
\begin{equation*}
Z_{(P,P')}\subset\{ \rho\in Z_{P}\ | \ 
\text{$\tilde q_j(\rho)=0 \iff n_j\in I^{P'}$} \},
\end{equation*}
and the two sets are not equal.  
In particular 
$Z_P\ne \bigsqcup_{P'\supseteq P} Z_{(P,P')}$, see
for example Remark~\ref{r:noninj}.

\section{Mirror conjecture and quantum cohomology}\label{s:MS}
The Givental/Eguchi-Hori-Xiong type mirror conjecture \cite[Conjecture 5.5.1]{BCKS:MSPFlag}
associated to the data 
introduced in the previous section
states that a full set of solutions to the quantum cohomology
$D$-module (e.g. \cite[Chapter~10]{CoxKatz:QCohBook}) 
associated to $SL_{n+1}/P$ can be written down 
on the mirror side as 
complex oscillatory integrals of the form
 $$
S_\Gamma(s):=\int_{\Gamma_s} e^{\mathcal F/z}\omega_s.
 $$
Here $s=(s_1,\dotsc, s_{k+1})\in \C^{k+1}$. Furthermore
$\omega_s$ is
a particular volume form on
$Z_{\tilde Q(s)}:=Z_{(e^{s_1-s_2},
\dotsc, e^{s_k-s_{k+1}})}$, and 
$\Gamma$ is a suitable family of (possibly
non-compact) middle-dimensional cycles
$\Gamma_s\subset Z_{\tilde Q(s)}$ for which the integral 
converges. In the $SL_{n+1}/B$ 
case this conjecture was proved by Givental~\cite{Giv:MSFlag}.

Whenever the conjecture holds 
the variety swept out by the critical points of 
$\mathcal F$ along the fibers of $Z^\circ\to (\C^*)^k$ should 
satisfy the relations of the small quantum cohomology ring 
(compare  with \cite{Giv:EquivGW}), or ideally 
completely recover the spectrum 
$Spec\left (qH^*(SL_{n+1}/P)[q_1\inv,\dots, q_k\inv]\right)$.

\subsection{}\label{s:crit}
Let 
$$Z^{\circ,crit}=Z^{\circ,crit}_P:=
\left\{\rho\in Z^\circ\ | \text{ $ \mathcal F|_{Z_{\tilde q(\rho)}}$ has
a critical point at $\rho$}  \ \right\}.$$ 

Following \cite{Giv:MSFlag} we write $\mathcal F$ in 
logarithmic vertex variables $T_v\in\C$ with $e^{T_v}=t_v$ to obtain
$$
\frac{\partial}{\partial T_v}\mathcal F=
\sum_{a, h_a=v} e^{T_{h(a)}-T_{t(a)}}-
\sum_{a, t_a=v} e^{T_{h(a)}-T_{t(a)}}.
$$ 
Therefore the critical point condition reads 
\begin{equation}\label{e:crit}
\sum_{a, h(a)=v} a-\sum_{a, t(a)=v} a=0, \qquad \text{ for 
all $v\in\mathcal V_\bullet$.} 
\end{equation}
For every vertex in $\mathcal V_\bullet$, the sum of incoming 
variables equals the sum of outgoing variables. 
We define 
\begin{equation}
Z^{crit}=Z^{crit}_P:=\left\{\rho\in Z_P\ |\ \text{ $\rho$ satisfies 
\eqref{e:crit}}\right\}. 
\end{equation}

\section{The GBCKS construction and the Peterson variety}\label{s:main}

In this section we demonstrate explicitly how the GBCKS 
construction
relates to $qH^*(SL_{n+1}/P)$. This is best done by 
comparing $Z^{crit}$ with the Peterson variety $Y_P$. While 
in the full flag variety case $Z^{\circ,crit}$
is almost isomorphic to $Y_B^\circ$, or $Spec(qH^*(SL_{n+1}/B)[q_1\inv,
\dotsc,q_n\inv])$ (it is isomorphic to an 
open dense subset), we will see 
that in the partial flag variety case
entire irreducible components of $Y_P^\circ$ can be missed out
by $Z^{\circ,crit}$. Nevertheless our result, see in particular
Proposition~\ref{p:qcoh}, should be considered as positive 
evidence for the mirror conjecture from \cite{BCKS:MSPFlag}. 
Although, as it turns out, the GBCKS mirror family $Z$, or rather $Z^\circ$, 
may be thought of as an open subset of a more 
complete (and canonical) mirror family, see \cite{Rie:MSgen}.

\subsection{}\label{s:g} We want to define a map $\phi:Z\to SL_{n+1}/B^-$. Let us 
first introduce some new notation. Set
\begin{equation}
l_j:= n_j-n_{j-1}, \quad\text { for $j=1,\dotsc, k+1$.}
\end{equation}
For $(m,r)\in \mathcal V_{\bullet}$ let
\begin{equation}
\tilde c_{m,r}=\begin{cases} c_{n_j, n_{j-1}+1}\prod^{p}_{i=2}
d_{n_j,n_{j-1}+i} & \text{if $m=n_j$ and $r=n_{j-1}+p$,}\\
&\text{where $2\le p\le l_j$ and $j=1,\dotsc,k$,}\\
c_{m,r} & \text{otherwise.}
\end{cases}
\end{equation}
Note that $\deg(\tilde c_{n_j,n_{j-1}+p})=p$. We also define for later use
\begin{equation}
\tilde d_{m,r+1}=\begin{cases} d_{n_j-1,n_{j-1}+1}
\prod_{i=1}^{p-1} c_{n_j-i, n_{j-1}} & \text{if $r=n_{j-1}$ and 
$m=n_j-p$,}\\
&\text{where $2\le p\le l_j$ and $j=1,\dotsc, k$,}\\ 
d_{m,r+1} & \text{otherwise.}  
\end{cases}
\end{equation}
Note that $\tilde q_j=\tilde c_{n_j, n_{j}}\tilde d_{n_j, n_j+1}$.  

Consider the simple root subgroups
$x_i(t)$ for $i\in I$. 
Let us also fix a one-parameter 
subgroup associated to a positive root 
$\alpha_{[i,i']}:=\alpha_i+\alpha_{i+1}+\dotsc+\alpha_{i'}$ with 
$1\le i<i'\le n$ by defining
\begin{equation}\label{e:1PS}
x_{[i,i']}(t):=\dot s_{i'}\dot s_{i'-1}\dotsc \dot s_{i+1}x_{i}(t)
\dot s_{i+1}\inv\dotsc \dot s_{i'-1}\inv\dot s_{i'}\inv.
\end{equation}
Explicitly,  $x_{[i,i']}(t)$ is the unipotent upper-triangular matrix with 
$(i, i'+1)$-entry $t$ and zeros everywhere else above the diagonal.

Let  $r=n_{j-1}+p $ for some $j=1,\dotsc, k$ and $1\le p\le l_j$. 
We define elements $g_r$ in $SL_{n+1}(\C[Z])$ by
\begin{eqnarray*}
g_{n_{j-1}+1}&=&x_{n}(c_{n,n_{j-1}+1})x_{n-1}(c_{n-1,n_{j-1}+1})
\dotsc x_{n_j}(c_{n_j,n_{j-1}+1})\dot s_{n_{j}-1}\dotsc\dot s_{n_{j-1}+1},\\
&\vdots&\\
g_r&=&
x_{n}(c_{n,r}) x_{n-1}(c_{n-1, r})\dotsc\dotsc 
x_{n_{j}}(\tilde c_{n_j,r})
\dot s_{n_{j}-1}\dot s_{n_{j}-2}\dotsc \dot s_{n_{j-1}+p},\\
&\vdots&\\
g_{n_j}&=&x_{n}(c_{n,n_j})x_{n-1}(c_{n-1,n_j})
\dotsc x_{n_j}(\tilde c_{n_j+1, n_j}).
\end{eqnarray*}
The element $g_r$ should be viewed as associated to the $r$-th column in the 
graph $(\mathcal V,\mathcal A)$.
For  $r= n_k+p$ with $1\le p\le l_{k+1}-1$ set 
\begin{eqnarray*}
g_{n_k+p}=\dot s_n \dot s_{n-1}\dotsc \dot s_{n_k+p}.
\end{eqnarray*} 
We can now form the product to get a new element
$g:=g_1 g_2\dots g_n \in SL_{n+1}(\C[Z])$, 
or equivalently a map
\begin{equation*}
g:Z\to SL_{n+1}(\C).
\end{equation*}
We define the map $\phi$, or $\phi_P$, keeping track of the dependence 
on $P$, by
\begin{equation}
\label{e:phi}
\begin{array}{llcl}
\phi_P: & Z_P&\to& SL_{n+1}/B^-,\\
&\rho&\mapsto& g(\rho)B^- .
\end{array}
\end{equation}
Note that the image of  $\phi_P$  
lies in $B^+\dot w_P B^-/B^-$. 

\subsection{Deodhar strata}\label{s:Deodhar} 
The intersections of opposite Bruhat cells $\mathcal R_{v,w}$ have
a decomposition into finitely many strata (each of the form $\C^l\x(\C^*)^m$)
due to Deodhar \cite{Deo:decomp}. We will not give Deodhar's original 
definition here, but rather use an equivalent description from 
\cite{MarRie:ansatz} which is ideally suited to our needs. 

The Deodhar decomposition 
of $\mathcal R_{v,w}$ depends on a choice of reduced 
expression for the longer element, $w$. We write 
$\mathbf w=s_{i_1}\dotsc s_{i_m}$ to mean $w$ with the given 
reduced expression $(i_1,\dotsc, i_m)$. A  
sequence of integers $1\le j_1<\dotsc<j_l\le m$ gives rise to
a subexpression $\mathbf v$ for $v$ in $\mathbf w$ if 
$s_{i_{j_1}}s_{i_{j_2}}\dotsc s_{i_{j_l}}=v$. The latter need not be 
a reduced expression for $v$. Associated to the pair $(\mathbf v,\mathbf w)$
of reduced expression $\mathbf w$ and subexpression $\mathbf v$ 
we have the sets
 \begin{eqnarray*}
J^+_{(\mathbf v,\mathbf w)}&=& 
\left\{ r=j_p\ | \ \text{some $p=1,\dotsc, l$ with 
$s_{j_1}\dotsc s_{j_{p-1}}>s_{j_1}\dotsc s_{j_{p-1}}s_{j_p}$ }
\right\},\\
J^-_{(\mathbf v,\mathbf w)}&=&\left\{ r=j_p\ | \ 
\text{some $p=1,\dotsc, l$ with 
$s_{j_1}\dotsc s_{j_{p-1}}<s_{j_1}\dotsc s_{j_{p-1}}s_{j_p}$ }
\right\},\\
J^\circ_{(\mathbf v,\mathbf w)} &=&\{1,\dotsc, m\}\setminus\{j_1,\dotsc, j_l\}.
 \end{eqnarray*}

The strata of $\mathcal R_{v,w}$ are indexed
by certain subexpressions $\mathbf v$ for $v$ 
in $\mathbf w$ called distinguished,
see \cite{Deo:decomp} or \cite[Section~3]{MarRie:ansatz}
for a definition. By \cite[Proposition~5.2]{MarRie:ansatz} 
the Deodhar stratum corresponding to 
$\mathbf v,\mathbf w$ is given by
\begin{equation}\label{e:DeoComp}
\mathcal R_{\mathbf v,\mathbf w}=
\left\{g_1 g_2\dotsc g_mB_-\ \left | \ g_r=\begin{cases}
\dot s_{i_r} & \text{if $r\in J^+_{(\mathbf v,\mathbf w)}$,}\\
y_{i_r}(m_r)\dot s_{i_r}\inv,\ \, 
m_r\in\C, &\text{if $r\in J^-_{(\mathbf v,\mathbf w)}$},\\
x_{i_r}(t_r),\qquad\quad t_r\in\C^*, &
\text{if $r\in J^\circ_{(\mathbf v,\mathbf w)}$}
\end{cases}\quad
\right.\right\}.
\end{equation}
Moreover the parameters $t_r\in\C^*$ and $m_r\in\C$ define an
isomorphism 
$$
(\C^*)^{J^\circ_{(\mathbf v,\mathbf w)}}\x
\C^{J^-_{(\mathbf v,\mathbf w)}}
\overset\sim\To\mathcal R_{\mathbf v,\mathbf w}.
$$
There is a unique distinguished subexpression $\mathbf v^+$ of $\mathbf w$
with $J^-_{(\mathbf v^+,\mathbf w)}=\emptyset$, which we call the positive
subexpression for $v$ in $\mathbf w$. It can be constructed as the rightmost 
reduced subexpression for $v$ in $\mathbf w$, see for example 
\cite[Lemma~3.5]{MarRie:ansatz}, and it corresponds to the 
unique open stratum $\mathcal R_{\mathbf v^+,\mathbf w}$
in $\mathcal R_{v,w}$.

\subsection{}
Consider the reduced expression $\mathbf {w_0}$ of $w_0$ given by
$$
(s_{n} s_{n-1}\dotsc s_1)(s_n s_{n-1}\dotsc s_2)\cdots (s_n s_{n-1}) s_n. 
$$
Let $P'$ be a parabolic with $0\le a_1<b_1<a_2<b_2<\dotsc<
a_h <b_h\le n$ such that
$$
I_{P'}=[a_1+1,b_1]\cup [a_2+1,b_2]\cup\cdots\cup [a_{h}+1, b_h], 
$$ 
as union of intervals in $\{1,\dotsc,n\}$. We 
have a reduced expression $\mathbf {w_{P'}}$ given by  
\begin{align*}
&(s_{b_1}s_{b_{1}-1}\dotsc s_{a_{1}+1})
(s_{b_1}s_{b_{1}-1}\dotsc s_{a_{1}+2})
\cdots (s_{b_{1}}s_{b_{1}-1}) s_{b_1}\\
&(s_{b_2}s_{b_2-1}\dotsc s_{a_2+1})
(s_{b_2}s_{b_2-1}\dotsc s_{a_2+2})
\cdots (s_{b_2}s_{b_2-1}) s_{b_2}\\
&\qquad\qquad\qquad\qquad\qquad\qquad\cdots\\
&(s_{b_h}s_{b_{h}-1}\dotsc s_{a_{h}+1})
(s_{b_h}s_{b_{h}-1}\dotsc s_{a_{h}+2})
\cdots (s_{b_{h}}s_{b_{h}-1}) s_{b_h}.
\end{align*}
The expression 
$\mathbf {w_{P'}}$ can also be constructed as the reduced expression
obtained from the positive subexpression 
for $w_{P'}$ in $\mathbf{w_0}$.

\begin{lemdef}\label{l:maps}
Let  $P'$ be a parabolic subgroup with $P'\supseteq P$ and recall
the definition of $Z_{(P,P')}$ from Section~\ref{s:Zparts}.
The map $\phi_P:Z_P\to G/B^-$ from \eqref{e:phi} restricts to 
$$
\phi_{(P,P')}: Z_{(P,P')}\To 
\mathcal R_{\mathbf{w_P^+},\mathbf{w_{P'}}}.
$$
In particular setting $P'=G$ gives $Z_{(P,G)}=Z^\circ_P$ and 
we define
$$
\phi^\circ_P:=\phi_{(P,G)}: Z^\circ_P\To
\mathcal R_{\mathbf{w_P^+},\mathbf{w_0}}.  
$$
\end{lemdef}

\begin{proof}[Proof of Lemma~\ref{l:maps}]
This lemma follows directly from the definitions of 
$Z_{(P,P')}$ and the map
$\phi_P$ together with the description of the Deodhar
strata proved in \cite[Proposition~5.2]{MarRie:ansatz}, 
see \eqref{e:DeoComp}. 
\end{proof}

We can now use $\phi_P$ to relate the GBCKS construction to 
the Peterson variety.

\begin{thm}\label{t:theorem1}
\begin{enumerate}\item
The map $\phi_P$ restricts to a map
\begin{equation*}
\phi_P^{crit}:Z^{crit}_P\To Y_P
\end{equation*} 
such that the following diagram commutes,
\begin{equation*}
\begin{matrix}
Z^{crit}_P&\overset{\phi^{crit}_P}\To  & Y_P\\
\ \underset{\tilde q\ \ }\searrow &   &  \underset{\ \ q}\swarrow\quad \\
& \C^{k} &              \qquad .
\end{matrix}
\end{equation*}
Here $q: Y_P\to \C^k$ is the map given by 
the quantum parameters $q_1,\dotsc, q_k$ in $qH^*(G/P)$ after applying
Peterson's isomorphism $\psi_P\inv$, see Theorem~\ref{t:Peterson}. 
\item
The morphism $\phi_P^{crit}$ restricted to the sets $Z^{crit}_{(P,P')}$ 
gives rise to embeddings
\begin{eqnarray*}
\phi_{(P,P')}^{crit}: Z^{crit}_{(P,P')}\To Y_{(P,P')}.
\end{eqnarray*}
The image of $\phi_{(P,P')}^{crit}$ is the 
intersection of $Y$ with the open Deodhar stratum 
$\mathcal R_{\mathbf{w_P^+},\mathbf{w_{P'}}}$ inside 
$\mathcal R_{w_P,w_{P'}}$, and we have an isomorphism
\begin{equation}\label{e:iso}
Z^{crit}_{(P,P')}\overset\sim\To Y\x_{G/B^-}
\mathcal R_{\mathbf{w_P^+},\mathbf{w_{P'}}}.
\end{equation}
\end{enumerate}
\end{thm}

\begin{rem}\label{r:noninj} 
The map $\phi_P^{crit}$ is not injective
outside the special subsets $Z_{(P,P')}$. 
For example for $SL_3/B$ consider the one-parameter 
family inside $Z^{crit}$ given by assigning values to the arrows
in $\mathcal A$ as follows
\begin{equation*}\rho_x=\quad
\begin{matrix}
\quad\star      & & & \\
\, 0 \uparrow & & & \\
\quad \bullet   &\overset{-x}\longrightarrow & \star\qquad\ \ & \\
x \uparrow &   & \uparrow -x\ \ \ & \\
\quad\bullet  & \overset{x}\longrightarrow & \bullet
\ \ \overset 0\longrightarrow &\star
\end{matrix}.
\end{equation*}
Then $\phi^{crit}_B(\rho_x)=B^-$ for all $x\in\C$. 
Note that $\rho_x$ does not lie in 
$\bigsqcup_{P'}Z_{(B,P')}$ unless $x=0$. 
\end{rem}

\section{Proof of Theorem~\ref{t:theorem1}}\label{s:proof1}

\subsection{}\label{s:recursions}
To prepare for proving the theorem
we first require some more notation and a 
technical lemma. We have fixed the parabolic $P$. 
Let $\mathcal I^{crit}$ denote the ideal  in $\C[Z]$ 
generated by the critical point conditions \eqref{e:crit}.
We set $d_{m,r}=0$ if $(m,r)\notin\mathcal V$
or $r=1$.

Let $(m,r)\in \mathcal V_\bullet$ and $l\ge 0$. 
Then to any set of columns $1\le r_1<r_2<\dotsc
<r_s\le r$ associate rows $m_1>m_2>\dotsc >m_s$ by $m_s=m$, and 
$m_{i-1}=m_{i}-deg(\tilde c_{m_{i},r_{i}})$. With this in mind let 
\begin{equation*}
G^{(m,r)}_l=\sum_{\begin{array}{c}
1\le r_1<\dotsc< r_s\le r\\
\sum deg(\tilde c_{m_i,r_i})=l
\end{array}}\left(\prod_{i=1}^s \tilde c_{m_i,r_i}\right ), 
\end{equation*}
if $l>0$, and set $G^{(m,r)}_0=1$. 
If $(m,r)\in \Z^2$ is not in $\mathcal V_{\bullet}$ then
we set $G^{(m,r)}_l=0$ by default. Also $G^{(m,r)}_l=0$ unless $l\le r$.

Recall the definition of $g$ from Section~\ref{s:g}
and let $u:=g\dot w_P\inv\in SL_{n+1}(\C[Z])$.
The element $u$ lies in $U^+(\C[Z])$ and is given by $u=u_1 u_2\dotsc u_{n_k}$
where for $n_{j-1}< r=n_{j-1}+p\le n_j$ we set
\begin{equation}\label{e:ur}
u_r=u_{n_{j-1}+p}=x_{n}(c_{n,r}) x_{n-1}(c_{n-1, r})\dotsc 
x_{n_j+1}(c_{n_j+1,r})x_{[n_{j}-p+1,n_j]}(\tilde c_{n_j,r}),
\end{equation}
with $j=1,\dotsc, k$, see \eqref{e:1PS}. Multiplying together the factors
$u_1\dotsc u_{n_k}$ 
it is straightforward  to check that $u$ is  
the $(n+1)\x (n+1)$-matrix 
\begin{equation}\label{e:u}
u=\left (U^{(0)}| U^{(1)} |\cdots | U^{(k)}\right )
\end{equation}
where $U^{(j)}$ is the $(n+1)\x l_{j+1}$ matrix
given explicitly by
\begin{equation*}
U^{(j)}=
\begin{pmatrix}
G^{(n_{j}, n_{j})}_{n_j} & 0    &    &     &     \\
G^{(n_j,n_{j})}_{n_j-1} &G^{(n_j+1,n_j)}_{n_{j}} &&&\\
 \vdots  & G^{(n_j+1,n_{j})}_{n_j-1} 
                     &  &\ddots &0\\
 \vdots   & & &\ddots &  G^{(n_{j+1}-1,n_j)}_{n_j}\\
 G^{(n_j,n_j)}_1 &&&&G^{(n_{j+1}-1,n_j)}_{n_j-1} \\
 1 & G^{(n_{j}+1,n_j)}_1 &&&\vdots \\
0 & 1 &       &&           \vdots  \\
  &   &       &        & G^{(n_{j+1}-1,n_j)}_{1}  \\
  &   &       &        &        \\
&&&&1\\
&&&& \\
&&&&\\
&&&&\\
0 &&&&0
\end{pmatrix}.
\end{equation*}
Note that $U^{(0)}$ is zero above the 
diagonal.
In general $G^{(m,r)}_l$ is a
matrix entry in the partial product $u_{(r)}=u_1 u_2\dotsc u_r$.

The definition of $G^{(m,r)}_l$ implies the following recursion.
\begin{equation}\label{e:recursion1}
  G^{(m,r)}_{l}=G^{(m,r-1)}_l+ \tilde c_{m,r} G^{(m-p,r-p)}_{l-p}
\end{equation}
where $p:=deg(\tilde c_{m,r})$. 

\begin{lem}   
If $(m,r)\in \mathcal V_{\bullet}$ and 
$0\le l\le r$ then 
\begin{equation}\label{e:recursion2}
G_l^{(m,r)}=G_l^{(m+1,r)}+d_{m,r+1} G_{l-1}^{(m,r-1)} \mod \mathcal I^{crit}. 
\end{equation}
\end{lem}

\begin{proof}[Proof of the Lemma] If $r=1$ then $l=0,1$ and 
the relation~\eqref{e:recursion2} is either trivial or it reads 
$c_{m,1}= c_{m+1,1}+ d_{m,2} $,
which is precisely the critical point condition at the vertex 
$(m,1)$. We now proceed by induction on $r$. The 
equalities in this proof are meant modulo $\mathcal I^{crit}$. 

We apply the induction hypothesis to the summands on the right 
hand side of \eqref{e:recursion1} to obtain
\begin{multline}\label{e:induction}
G^{(m,r)}_l=G^{(m+1,r-1)}_{l}+d_{m,r} G^{(m,r-2)}_{l-1}
+ \tilde c_{m,r} G^{(m-p+1,r-p)}_{l-p}\\
+\tilde c_{m,r } d_{m-p,r-p+1} G_{l-p-1}^{(m-p,r-p-1)},
\end{multline}
where $p$ is fixed to be the degree of $\tilde c_{m,r}$. 
\vskip 3mm
\noindent{\it Case 1~:} Suppose $\tilde c_{m,r}=c_{m,r}$. 
Then we can substitute
$$
c_{m,r}=c_{m+1,r}-d_{m,r}+d_{m,r+1}\qquad
\text{and} \qquad
c_{m,r}d_{m-1,r}= d_{m,r}c_{m,r-1}
$$
to obtain
\begin{eqnarray*}
G^{(m,r)}_{l}=G^{(m+1,r-1)}_l+c_{m+1,r}G^{(m,r-1)}_{l-1}
+ d_{m, r+1}G^{(m,r-1)}_{l-1}+ \\
+d_{m,r}\left(G^{(m,r-2)}_{l-1}-G^{(m,r-1)}_{l-1}+
c_{m,r-1} G^{(m-1,r-2)}_{l-2} \right).
\end{eqnarray*}
Now $\tilde c_{m,r}=c_{m,r}$ implies also 
$\tilde c_{m,r-1}=c_{m, r-1}$ and
$\tilde c_{m+1,r}=c_{m+1,r}$. Therefore \eqref{e:recursion1}
applies twice to give 
\begin{equation*}
G^{(m,r)}_{l}=G_l^{(m+1,r)}+d_{m,r+1} G^{(m,r-1)}_{l-1}.
\end{equation*}
\noindent{\it Case 2~:} Suppose $(m,r)=(n_j,n_{j-1}+p)$ for
some $j=2,\dotsc,k$ and $1\le p\le l_j$. 
In this case the vertex 
$(m-p,r-p)$ lies on the right hand edge of 
the graph, and
$d_{m-p,r-p+1}=0$. Furthermore by the critical point condition
at the vertex $(m,r)$ we can substitute $d_{m,r}=c_{m+1, r}+d_{m,r+1}$.
So \eqref{e:induction} becomes
\begin{multline*}
G^{(m,r)}_l=G^{(m+1,r-1)}_{l}+(c_{m+1, r}+d_{m,r+1}) G^{(m, r-2)}_{l-1}
+\tilde c_{m,r}G^{(m-p+1,r-p)}_{l-p}\\
=G^{(m+1,r-1)}_l+ c_{m+1,r}\left (G^{(m,r-1)}_{l-1}-\tilde c_{m,r-1}
G^{(m-p+1,r-p)}_{l-p}\right) + d_{m,r+1}G^{(m, r-2)}_{l-1}
\\ +
\tilde c_{m,r}G^{(m-p+1,r-p)}_{l-p}\\
=G^{(m+1,r)}_l- c_{m+1,r}\tilde c_{m,r-1}G^{(m-p+1,r-p)}_{l-p} + d_{m,r+1} 
G^{(m,r-2)}_{l-1}+\tilde c_{m,r} G^{(m-p+1,r-p)}_{l-p}.
\end{multline*}
Finally we substitute $\tilde c_{m,r}=\tilde c_{m, r-1} d_{m,r}=
\tilde c_{m, r-1}(c_{m+1,r}+ d_{m, r+1})$ to get
\begin{multline*}
G^{(m,r)}_l=G^{(m+1,r)}_l+\tilde c_{m,r-1}d_{m,r+1}G^{(m-p+1,r-p)}_{l-p}+
d_{m,r+1} G^{(m,r-2)}_{l-1}\\
= G^{(m+1,r)}_l+
d_{m,r+1}\left(G^{(m,r-2)}_{l-1}+\tilde c_{m, r-1}G^{(m-p+1,r-p)}_{l-p}
\right)\\
=G^{(m+1,r)}_{l}+ d_{m,r+1} G^{(m,r-1)}_{l-1}.
\end{multline*}
\end{proof}

\begin{cor}\label{c:lem}
The elements $G^{(n_{j}+p,n_j)}_l$ appearing as matrix entries
in $u$ satisfy
\begin{equation}
\label{e:shaperel}
G^{(n_j,n_j)}_l=G^{(n_j+1,n_j)}_l=\cdots=G^{(n_{j+1}-1,n_j)}_l \mod \mathcal I^{crit}.
\end{equation}
In other words the $U^{(j)}$  as matrices of functions on $Z^{crit}$ are 
constant along the diagonals. \qed
\end{cor}

\subsection{}
We now use the results from Section~\ref{s:recursions}
to show that the elements $G^{(n_j,n_j)}_{l}|_{Z^{crit}}$
in $\C[Z^{crit}]$ satisfy the
relations of the $E^{(j)}_l$ in $qH^*(SL_{n+1}/P)$. 
\begin{prop}\label{p:qcoh}
The assignments 
\begin{align*}
q^P_j&\mapsto \tilde c_{n_j,n_j}\tilde d_{n_j, n_j+1}&  & \text{for $j=1\dotsc, k$, and}\\
\sigma^{(j)}_p&\mapsto 
\tilde c_{n_j,n_{j-1}+p}+(-1)^{p}\tilde d_{n_j-p, n_{j-1}+1}& &\text{for $j=1,\dotsc, k+1$,}
\end{align*}
where $1\le p\le l_j$, define an algebra homomorphism 
$\kappa:qH^*(SL_{n+1}/P)\To \C[Z^{crit}]$.
The homomorphism $\kappa$ takes $E^{(j)}_l$ to $G^{(n_j, n_j)}_{l}|_{Z^{crit}}$.
\end{prop}

\begin{proof}
Let 
$\tilde\sigma^{(j)}_{p}
:=\tilde c_{n_j,n_{j-1}+p}+(-1)^{p}\tilde d_{n_j-p, n_{j-1}+1}$
for $1\le p\le l_{j}$, and let $\tilde\sigma^{(j)}_0=1$. In all other
cases set $\tilde\sigma^{(j)}_{p}=0$.
It suffices to prove the relation
\begin{multline}\label{e:relation}
G^{(n_{j-1},n_{j-1})}_l= G^{(n_j,n_j)}_{l}-\left(
\tilde\sigma^{(j)}_1 G^{(n_j,n_j)}_{l-1} + 
\tilde\sigma^{(j)}_2 G^{(n_j,n_j)}_{l-2}+
\dotsc + \tilde\sigma^{(j)}_{l}
\right)\\ + (-1)^{l_j}q_{j-1} G^{(n_{j-2},n_{j-2})}_{l-n_j+n_{j-2}},
\end{multline}
$\mod \mathcal I^{crit}$, where $j=1,\dotsc, k+1$. See Section~\ref{s:qcoh}. 

Using Corollary~\ref{c:lem} we replace the left hand side of 
\eqref{e:relation} by $G^{(n_{j}-1,n_{j-1})}_l$ and then 
apply \eqref{e:recursion2} to get
\begin{equation}\label{e:twoterm}
G^{(n_{j-1},n_{j-1})}_{l}=
G^{(n_j,n_{j-1})}_l+d_{n_{j}-1,n_{j-1}+1} 
G^{(n_{j}-1,n_{j-1}-1)}_{l-1} \ \mod \mathcal I^{crit}.
\end{equation}
Now we 
consider the first summand and successively 
apply the relation \eqref{e:recursion1} 
\begin{multline}\label{e:first}
G^{(n_j,n_{j-1})}_l=
G^{(n_j,n_{j-1}+1)}_l-c_{n_j,n_{j-1}+1}G^{(n_j-1,n_{j-1})}_{l-1}
\\
=G^{(n_j,n_{j-1}+2)}_{l}-\tilde c_{n_j,n_{j-1}+2}
G^{(n_j-2,n_{j-1})}_{l-2} - 
c_{n_j,n_{j-1}+1}G^{(n_j-1,n_{j-1})}_{l-1}= \dotsc\\
\dotsc =G^{(n_j,n_j)}_l-
\sum_{i=1}^{n_j-n_{j-1}}\tilde c_{n_j,n_{j-1}+i} G^{(n_j-i,n_{j-1})}_{l-i}.
\end{multline} 
Let us apply the same relation to the second summand in 
\eqref{e:twoterm},
\begin{equation*}
d_{n_j-1,n_{j-1}+1}G^{(n_j-1,n_{j-1}-1)}_{l-1}=
d_{n_j-1,n_{j-1}+1}\left(G^{(n_j-1,n_{j-1})}_{l-1}-c_{n_j-1,
n_{j-1}}G^{(n_j-2,n_{j-1}-1)}_{l-2}\right),
\end{equation*}
and note that we can make the replacement 
$d_{n_{j}-1,n_{j-1}+1}c_{n_j-1,n_{j-1}}=\tilde d_{n_j-2,n_{j-1}+1}$.

Repeating this process, successively applying \eqref{e:recursion1}
to the final summand we get
\begin{multline}\label{e:second}
d_{n_{j-1},n_{j-1}+1}G^{(n_{j}-1, n_{j-1}-1)}_{l-1}\\
\ \ =d_{n_{j}-1,n_{j-1}+1}G^{(n_j-1,n_{j-1})}_{l-1}-\tilde d_{n_j-2,n_{j-1}+1}
G^{(n_{j}-2,n_{j-1}-1)}_{l-2}\\
=d_{n_{j}-1,n_{j-1}+1} G^{(n_j-1,n_{j-1})}_{l-1}-
\tilde d_{n_{j}-2,n_{j-1}+1} G^{(n_j-2,n_{j-1})}_{l-2} \\ + 
\tilde d_{n_{j}-3, n_{j-1}+1} G^{(n_j-3,n_{j-1}-1)}_{l-3}
=\dotsc \\
=d_{n_j-1,n_{j-1}+1}G^{(n_{j}-1,n_{j-1})}_{l-1}-\cdots \qquad\qquad\qquad\qquad\qquad\\
\qquad \qquad + (-1)^{n_j-n_{j-1}-1}\tilde d_{n_{j-1},n_{j-1}+1}
G^{(n_{j-1},n_{j-1}+1)}_{l-n_j+n_{j-1}}\\
= \left( \sum_{i=1}^{n_j-n_{j-1}}(-1)^{i+1}\tilde d_{n_j-i,n_{j-1}+1}
G^{(n_j-i,n_{j-1})}_{l-i}\right)\qquad \qquad\\ 
+ (-1)^{n_j-n_{j-1}}
\tilde d_{n_{j-1},n_{j-1}+1}\tilde c_{n_{j-1},n_{j-1}}
G^{(n_{j-2},n_{j-2})}_{l-n_{j}+n_{j-2}}.
\end{multline}
Summing \eqref{e:first} and \eqref{e:second} gives
\begin{multline*}
G^{(n_{j-1},n_{j-1})}_l=G^{(n_j,n_{j-1})}_l+d_{n_{j}-1,n_{j-1}+1} 
G^{(n_{j}-1,n_{j-1}-1)}_{l-1}\\
=G^{(n_j,n_j)}_{l}-\sum_{i=1}^{n_j-n_{j-1}}
\tilde\sigma^{(j)}_i G^{(n_j-i,n_{j-1})}_{l-i}+
(-1)^{n_j-n_{j-1}}
q_{j-1} G^{({n_{j-2}, n_{j-2})}}_{l-n_j+n_{j-2}}\ \mod\mathcal I^{crit}.
\end{multline*}
Using Corollary~\ref{c:lem} we see that this is 
the relation \eqref{e:relation} we were trying to prove. 
\end{proof}

\begin{rem}
Note that \eqref{e:first} and \eqref{e:second} were obtained 
using only the definition of the $G^{(m,r)}_l$. We see therefore that the following
relation, 
\begin{multline}\label{e:inZP}
G^{(n_{j+1},n_{j})}_l+d_{n_{j+1}-1,n_{j}+1} 
G^{(n_{j+1}-1,n_{j}-1)}_{l-1}\\
=G^{(n_{j+1},n_{j+1})}_l-\sum_{i=1}^{n_{j+1}-n_{j}}\left(\tilde c_{n_{j+1},n_{j}+i}
+(-1)^{i}\tilde d_{n_{j+1}-i,n_{j}+1}\right)
G^{(n_{j+1}-i,n_{j})}_{l-i}\\ 
+ (-1)^{n_{j+1}-n_{j}}
\tilde d_{n_{j},n_{j}+1}\tilde c_{n_{j},n_{j}}
G^{(n_{j-1},n_{j-1})}_{l-n_{j+1}+n_{j-1}},
\end{multline}
which is obtained by combining \eqref{e:first} and \eqref{e:second} and replacing 
$j$ by $j+1$, holds in $\C[Z_P]$. 
If $l>n_{j}$ then the left hand side of \eqref{e:inZP} is zero.
\end{rem}

We may now use these results to prove the theorem. 
For a different more Lie theoretic proof in the $G/B$
case see also \cite{Rie:MSgen}.

\begin{proof}[Proof of Theorem~\ref{t:theorem1}]
Consider the matrix $u\in U^+(\C[Z])$ from \eqref{e:u}
and let $\rho\in Z^{crit}$. A direct calculation using
the shape of $u$ (see Corollary~\ref{c:lem}) 
and the relation \eqref{e:relation}
shows that $u(\rho)\inv\cdot f\in
\dot w_P\cdot(\mathfrak b^-\oplus \sum_{i\in I}\C e_i)$
as required (compare \cite[Section~4.2]{Rie:QCohPFl}). So 
we have $\phi_P^{crit}: Z^{crit}\to Y_P$. 

Next we can evaluate the function from \eqref{e:M}
at $\phi_P(\rho)$ to get 
$$
M_{w_l^{[n_j]}\omega_{n_j}}(\phi_P(\rho))=
\left<u(\rho) \cdot v_{-\omega_{n_j}},\dot w^{[n_j]}_l\cdot 
v_{-\omega_{n_j}}\right>,
$$
for $j=1,\dotsc, k$. It follows from this that
$M_{w_l^{[n_j]}\omega_{n_j}}(\phi_P(\rho))=G^{(n_j,n_j)}_{l}(\rho)$.
 Therefore the map $(\phi_P^{crit})^*:\C[Y_P]\to \C[Z^{crit}]$ 
is the composition of Peterson's isomorphism 
$\psi_P:\C[Y_P]\to qH^*(SL_{n+1}/P)$ 
with the homomorphism $\kappa:qH^*(SL_{n+1}/P)\to\C[Z^{crit}]$ from 
Proposition~\ref{p:qcoh}. Since $\kappa$ takes $q_j$ to $\tilde q_j|_{Z^{crit}}$
this implies also the second part of $(1)$.

Let $I^{P'}=\{n_{j_1},\dotsc, n_{j_{k'}}\}\subset I^P$ and set $j_0=0$
and $j_{k'+1}=k+1$.  The variety $Z_{(P,P')}$ is 
isomorphic to a product of varieties 
$Z_{(P_i, SL_{l'_i})}$, where $l'_i=n_{j_{i}}-n_{j_{i-1}}$ for $i=1,\dots, k'+1$,
and the parabolic $P_i$
in $SL_{l'_i}$ is determined by $I^{P_i}=\{n_{j_{i-1}+1}-n_{j_{i-1}}, n_{j_{i-1}+2}-n_{j_{i-1}},\dotsc, n_{j_{i}-1}-n_{j_{i-1}}\}$.
On the other hand
we have corresponding coordinate projections 
$$
\mathcal R_{\mathbf {w_P^+},\mathbf {w_{P'}}}
\to \mathcal R^{SL_{l'_i}}_{\mathbf {w_{P_i}^+},\mathbf {w_0}},
$$ 
which are easily seen to be compatible with intersecting with the Peterson variety
(of $SL_{n+1}$ and $SL_{l'_i}$, respectively).  
In this way the problem of finding an 
inverse to $\phi^{crit}_{(P,P')}$ is reduced to finding inverses to the maps $\phi^{crit}_{(P_i,SL_{l'_i})}$.
Therefore we assume from now on that $P'=G$.

Suppose $\rho\in Z_{(P,G)}$ and 
$\phi_P(\rho)=g(\rho) B^-$ lies in $Y_{(P,G)}$. We can
recover the values $c_{i,j}(\rho)$ for all the vertical arrows 
from the factors of $\bar g:=g(\rho)$. 
(Recall that the entries of the simple root subgroup factors
in $\bar g$ 
are coordinates on the
Deodhar stratum where $\bar gB^-$ lies). From the special entries
$\tilde c_{n_j, n_{j-1}+p}(\rho)$ we also recover the 
values of particular horizontal arrows 
from the rim of the graph, namely the $d_{n_j,n_{j-1}+p}(\rho)$. 
Finally, from $q_j(\bar gB^-)$ (along with all of the
other coordinates already determined) we can work out 
values for the remaining horizontal arrows from the rim, 
the $d_{\star_{j+1}}(\rho)$. By the box equations 
\eqref{e:box} these values for all of the vertical 
arrows and for the rim determine a unique element in $Z_{(P,G)}$.
If $\rho$ was in $Z_{(P,G)}^{crit}$, then this element is precisely $\rho$.  
Applying the same procedure to an arbitrary element of 
$Y_P\x_{G/B^-} \mathcal R_{\mathbf {w_P^+},\mathbf{w_{G}}}$
defines a morphism 
$$
\beta:Y_P\x_{G/B^-} \mathcal R_{\mathbf {w_P^+},\mathbf{w_{G}}}
\To Z_{(P,G)},
$$ 
such that $\beta\circ \phi^{crit}_{(P,G)}$ 
is the identity on $Z^{crit}_{(P,G)}$. 
It remains to show that the image of $\beta$ lies in $Z^{crit}_{(P,G)}$.
Then $\beta$ is the inverse to \eqref{e:iso}
and (2) follows.

Consider  $\rho\in Z_{(P,G)}$ in the image of $\beta$. So $u(\rho)\dot w_P B^-\in Y_P$
and $\rho=\beta(u(\rho)\dot w_P B^-)$. Therefore $\rho$ satisfies an 
identity of $(n+1)\x (n+1)$ matrices over $\C[Z_{(P,G)}]$ of the 
following form,
\begin{equation}\label{e:PetersonCondition}
u (f+A_++ Q)=f u.
\end{equation}
Here $f$ is the principal nilpotent from \eqref{e:f} and $u$ the matrix from \eqref{e:u} with blocks 
$U^{(j)}$. 
The matrix $A_+$ is a block diagonal matrix
with upper-triangular blocks $A_+^{(j)}$ of size $l_j\x l_j$ for $j=0,\dotsc, k$, and $Q$ is the matrix
with entry $(-1)^{l_j}\tilde q_j$ in position $(n_{j-1}+1,n_j-1)$ for $j=1,\dotsc, k$ and zeroes elsewhere. 

We denote the $i$-th
column vector of $U^{(j)}$ by $U^{(j)}_i$. Let the entries of $A^{(j)}_+$ be denoted by $a_{r,m}^{(j)}$. 
The individual columns of \eqref{e:PetersonCondition} give identities
\begin{eqnarray}\label{e:columnm}
&U^{(j)}_{m+1}+ a_{1,m}^{(j)}U^{(j)}_1 + a_{2,m}^{(j)} U^{(j)}_2 + \dotsc + a_{m,m}^{(j)} U^{(j)}_m=f U^{(j)}_{m} ,
\\ \label{e:columnq}
&\quad U^{(j+1)}_{1}+ a_{1,l_j}^{(j)}U^{(j)}_1 +  \dotsc + 
a_{l_j,l_j}^{(j)} U^{(j)}_{l_j}+ (-1)^{l_j}\tilde q_j U^{(j-1)}_1=f U^{(j)}_{l_j},  
\end{eqnarray}
where $1\le m\le l_j-1$.

Note that 
\begin{equation}\label{e:Gmrr}
G^{(m,r)}_r=\prod_{i=1}^s \tilde c_{m_i,r_i},
\end{equation}
where $(m_1,r_1)=(m-r+1,1) $ and $(m_{i},r_{i})=(m_{i-1}+\deg(\tilde c_{m_i,r_i}),r_{i-1}+\deg(\tilde c_{m_i,r_i}))$. 
Therefore  $G^{(m,r)}_r$ is invertible in $\C[Z_{(P,G)}]$ and $G^{(m,r)}_r(\rho)\ne 0$, a fact we will use repeatedly without further mention. 

The identity \eqref{e:columnm}  implies
recursively that 
$$
a^{(j)}_{i,m}(\rho)=0\quad \text{and}\quad 
 G^{(n_j+{m-1},n_j)}_l(\rho)=G^{(n_j+m,n_j)}_l(\rho).
$$  
Similarly the identity \eqref{e:columnq} implies, 
that
$$
a^{(j)}_{i,l_j}(\rho) = -\tilde\sigma^{(j+1)}_{l_j-i+1}(\rho)
$$
for $1\le i\le l_j$, and 
$$
G^{(n_{j+1}-1,n_j)}_l(\rho)=G^{(n_{j+1},n_j)}_l(\rho)+d_{n_{j+1}-1,n_j+1}(\rho)G^{(n_{j+1}-1,n_j-1)}_{l-1}(\rho)
$$
for $1\le l\le n_j$ (comparing also with \eqref{e:inZP}). 
Therefore $\rho$ satisfies the relation \eqref{e:recursion2} at all vertices
$(m,r)=(n_j+p,n_j)$ with $0\le p\le l_j$. Moreover at the vertex $(n_j,n_j)$ and with $l=n_j$ this
relation reads
$$
G^{(n_j,n_j)}_{n_j}=G^{(n_j+1,n_j)}_{n_j},
$$
or equivalently,
$$
\tilde c_{n_1,n_1}\tilde c_{n_2,n_2}\dotsc \tilde c_{n_j,n_j}=
\tilde c_{n_1,n_1}\dotsc \tilde c_{n_{j-1},n_{j-1}} \tilde c_{n_j, n_j-1} c_{n_j+1,n_j}.
$$
Replacing $\tilde c_{n_j,n_j}$ by $\tilde c_{n_j,n_j-1}d_{n_j,n_j}$ and canceling we see 
therefore that $\rho$ satisfies 
\begin{equation*}
d_{n_j,n_j}(\rho)=c_{n_j+1,n_j}(\rho),
\end{equation*}
which is the critical point condition at ${(n_j,n_j)}$.

We will now prove using induction that $\rho$ satisfies the relation 
\eqref{e:recursion2} and the critical point condition for 
each of the remaining vertices in $\mathcal V_{\bullet}$.  
Let us consider the ordering on $\mathcal V_\bullet$ starting from $(n_1,n_1)$ and 
defined by
$(m',r') \le (m,r)$ if $m'<m$ or $m'=m$ and $r'\ge r$. We may assume 
that $\rho$ satisfies the relation \eqref{e:recursion2} and the critical point condition for all vertices 
$(m',r')$ and degrees $l'$ such that  
$(m',r')\le (m,r)$ and $l'<l$.

The start of induction at the vertex $(n_1,n_1)$ has already been checked. Let us prove the relation \eqref{e:recursion2} at a vertex $(m,r)$ which is not on the right hand edge of the graph (assuming as part of the induction hypothesis that everything is already proved for the right-most vertex in the row $m$). We have that  $(m,r+1)$ is another $\bullet$-vertex. Then at $\rho$,  
\begin{multline}
G^{(m,r)}_l=G^{(m,r+1)}_l-\tilde c_{m,r+1}G^{(m-p,r+1-p)}_{l-p}\\
=\left(G^{(m+1,r+1)}_l+d_{m,r+2} G^{(m,r)}_{l-1}
\right ) -\tilde c_{m,r+1}G^{(m-p,r+1-p)}_{l-p}\\
=\left(G^{(m+1,r)}_l+c_{m+1,r+1} G^{(m,r)}_{l-1} +d_{m,r+2} G^{(m,r)}_{l-1}\right )
 -\tilde c_{m,r+1}G^{(m-p,r+1-p)}_{l-p}\\
 = \begin{cases}
 G^{(m+1,r)}_l + (d_{m,r+1}+c_{m,r+1}) G^{(m,r)}_{l-1} - c_{m,r+1}G^{(m-1,r)}_{l-1}, & \text{if }
 \deg(\tilde c_{m,r+1})=1,\\
  G^{(m+1,r)}_l + d_{m,r+1} G^{(m,r)}_{l-1} -\tilde c_{m,r+1}G^{(m-p,r+1-p)}_{l-p},
  & \text{if }\deg(\tilde c_{m,r+1})=p>1.
 \end{cases}
\end{multline}
Here we used the induction hypothesis twice: first that \eqref{e:recursion2} holds and
then that the critical point condition holds at the vertex $(m,r+1)$. 

In the first of the two cases above we can go on to use the inductive assumption that \eqref{e:recursion2}
holds in degree $l-1$ at the vertex  $(m,r)$,  followed by a
box relation and \eqref{e:recursion1}, to obtain at $\rho$
\begin{multline}
G^{(m,r)}_l= G^{(m+1,r)}_l + d_{m,r+1} G^{(m,r)}_{l-1}+c_{m,r+1} (G^{(m,r)}_{l-1} - G^{(m-1,r)}_{l-1})\\
=G^{(m+1,r)}_l + d_{m,r+1} G^{(m,r)}_{l-1} - c_{m,r+1} d_{m-1,r+1} G^{(m-1,r-1)}_{l-2}\\
=G^{(m+1,r)}_l + d_{m,r+1} G^{(m,r)}_{l-1} - d_{m,r+1} c_{m,r} G^{(m-1,r-1)}_{l-2}\\
=G^{(m+1,r)}_l + d_{m,r+1} (G^{(m,r)}_{l-1} - c_{m,r} G^{(m-1,r-1)}_{l-2})=
G^{(m+1,r)}_l+d_{m,r+1} G^{(m,r-1)}_{l-1}.
\end{multline} 
 In the second case we have $ \tilde c_{m,r+1} =  d_{m,r+1} \tilde c_{m,r}$, and therefore again
 \begin{multline}
G^{(m,r)}_l=  G^{(m+1,r)}_l + d_{m,r+1} (G^{(m,r)}_{l-1} -\tilde c_{m,r}G^{(m-p,r+1-p)}_{l-p})
\\
= G^{(m+1,r)}_l+d_{m,r+1} G^{(m,r-1)}_{l-1}. 
 \end{multline}
So we see that $\rho$ satisfies  the relation $\eqref{e:recursion2}$ at the vertex $(m,r)$. 

Let us now show the critical point condition at $(m,r)$. Note first that since \eqref{e:recursion2} holds at 
$(m,r)$ we have
that $\rho$ satisfies
\begin{multline}\label{e:d-at-(m,r)}
G^{(m,r)}_r=G^{(m+1,r)}_r+d_{m,r+1}G^{(m,r-1)}_{r-1}\\
=c_{m+1,r}G^{(m,r-1)}_{r-1}+
d_{m,r+1}G^{(m,r-1)}_{r-1}=(c_{m+1,r}+d_{m,r+1})G^{(m,r-1)}_{r-1}
\end{multline}  
Now suppose first that  $(m,r)$ is of the form $(n_{j},n_{j}+s)$ and $2\le s\le l_{j+1}$.  
Then \eqref{e:Gmrr} implies that
$$
G^{(m,r)}_r=G^{(n_{j},n_{j}+s)}_{n_{j}+s}=d_{n_{j},n_{j}+s}G^{(n_{j},n_{j}+s-1)}_{n_{j}+s-1}.
$$
Comparing with \eqref{e:d-at-(m,r)} we see that 
$$
d_{n_{j},n_{j}+s}(\rho)= c_{n_{j}+1,n_{j}+s}(\rho)+d_{n_{j},n_{j}+s+1}(\rho),
$$
which is the critical point condition at $(n_{j},n_{j}+s)$. 

For all other vertices $(m,r)\in \mathcal V_\bullet$ we have $\tilde c_{m,r}=c_{m,r}$ and therefore
at $\rho$
\begin{multline}
G^{(m,r)}_r=c_{m,r}G^{(m-1,r-1)}_r=
c_{m,r}(G^{(m,r-1)}_{r-1}+d_{m-1,r}G^{(m-1,r-2)}_{r-2})\\
=c_{m,r}G^{(m,r-1)}_{r-1}+c_{m,r} d_{m-1,r}G^{(m-1,r-2)}_{r-2}=
c_{m,r}G^{(m,r-1)}_{r-1}+d_{m,r} c_{m,r-1}G^{(m-1,r-2)}_{r-2}\\
=(c_{m,r}+d_{m,r})G^{(m,r-1)}_{r-1}.
\end{multline}
Comparing this identity with \eqref{e:d-at-(m,r)} gives
$$
c_{m+1,r}(\rho)+d_{m,r+1}(\rho)=c_{m,r}(\rho)+d_{m,r}(\rho),
$$
which is the critical point condition at the vertex $(m,r)$.

The induction step showing the critical point condition also works for the 
vertices $(n_j+s,n_j)$ with $1\le s\le l_j$ along the 
right hand edge. Thus once we have proved the relations \eqref{e:recursion2} and the 
critical point conditions in the $m$-th row, where $m=n_j+s-1$, the critical point condition
at the vertex $(n_j+s,n_j)$ of the subsequent row follows, and allows us to continue the induction 
along $(m+1)$-st row.  Since we have already checked the critical point condition at all of the vertices
$(n_j,n_j)$, the induction now goes through to the end and 
implies that $\rho\in Z^{crit}_{(P,G)}$. In fact, by verifying the critical point conditions
for $\rho=\beta(gB^-)$ from the relations of the Peterson variety, we have shown that 
$\beta$ defines a morphism,  $Y_P\x_{G/B^-} \mathcal R_{\mathbf {w_P^+},\mathbf{w_{G}}}
\To Z^{crit}_{(P,G)}$  which is the desired inverse to \eqref{e:iso}. 
This completes the
proof.   
\end{proof}

\section{The example $Gr_2(\C^4)$.}\label{s:Gr2(4)}
\begin{figure}
\begin{center}
\leavevmode
 \[ \includegraphics{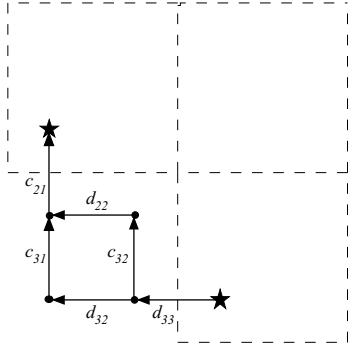}\]
\end{center}
\caption{The graph associated to $Gr_2(\mathbb C^4)$}\label{f:Gr24}
\end{figure}

Consider the mirror family for $Gr_2(\C^4)$ given in
\cite{BCKS:MSGrass,BCKS:MSPFlag}. 
It corresponds to the graph in Figure \ref{f:Gr24}.

In this case we have 
$$
\C[Z]=\C[c_{21},c_{31},c_{31}, d_{22},d_{32},d_{33}]/
(d_{32}c_{31}-c_{32}d_{22})
$$ 
and the critical point
condition (for $\tilde q= d_{33} c_{32} d_{22} c_{31}$ fixed) is
 $$
c_{21}=d_{22}+c_{31}, \qquad d_{33}=c_{32}+d_{32},\qquad c_{31}=d_{32},
\qquad c_{32}=d_{22}.
 $$

It is easy to check that the
critical point problem in this case has up to scalar only one solution:
 $$
c_{31}=d_{32}=c_{32}=d_{22}=1,\qquad c_{21}=d_{33}=2.
 $$
Or for fixed value of $q=1$, 
say, there are exactly $4$ solutions
 $$
c_{31}=d_{32}=c_{32}=d_{22}=
\exp(2\pi i k/4)/\sqrt{2},\qquad c_{21}=d_{33}=2\exp(2\pi i k/4)/\sqrt{2}
 $$
where $k=1,\dotsc,4$.
However $Spec(qH^*(Gr_{2}(\C^4)))$ over $q=1$  
should actually have $\dim H^*(Gr_2(\C^4))=6$ points.
(Note that by \cite{Gepner:FusRing} 
the quantum cohomology ring of a Grassmanian $Gr_d(\C^m)$
for fixed 
nonzero value of $q$ is semisimple. The 
$\left(d\atop m\right)$ points 
in the Peterson variety for any fixed value of $q$
are described in \cite{Rie:QCohGr}.)

We can find the two missing elements explicitly 
in the Peterson variety $Y_P$ for $Gr_{2}(\C^4)$. They are
$$
\begin{bmatrix}
1& 0 &i &0\\
 & 1 &0 &i\\
 &   &1 &0\\
 &   &  &1
\end{bmatrix}
\dot s_1\dot s_3 B^\vee_-/B^\vee_-, \quad\text{and}\quad 
\begin{bmatrix}
1& 0 &-i &0\\
 & 1 &0 &-i\\
 &   &1 &0\\
 &   &  &1
\end{bmatrix}
\dot s_1\dot s_3 B^\vee_-/B^\vee_-.
$$
The reason for the discrepancy is that $Y_P^\circ$ 
has two irreducible components,
one of them in 
the open Deodhar stratum $\mathcal R_{\mathbf {w_P^+},\mathbf w_0}$,
 $$
\left\{\begin{bmatrix}
1& 0 & s^2 & 0\\
 &  1  &\sqrt 2 s& s^2\\
 &     & 1     & \sqrt 2 s\\
 & & & 1 
\end{bmatrix}\dot s_1\dot s_3 B^\vee_-/B^\vee_- \ |\ s\in\C^*
\right\},
 $$
and thus captured by the GBCKS construction, the other, 
 $$
\left\{\begin{bmatrix}
1& 0 & m & 0\\
 &  1  &0 & m\\
 &     & 1     & 0\\
 & & & 1 
\end{bmatrix}\dot s_1\dot s_3 B^\vee_-/B^\vee_-\ |\ m\in\C^*
\right\},
 $$ 
in the smaller Deodhar stratum 
``$\,\mathcal R_{\bar 32\bar 1\bar 32\bar 3}\, $'' 
corresponding to the subexpression $s_3 1 s_1  s_3 1 s_3$ for $s_1 s_3$ in 
$s_3 s_2 s_1 s_3 s_2 s_3$. The elements in this stratum are of the form 
$$
\dot s_3 x_2(t_2) \dot s_1 y_3(m_4) \dot s_3\inv x_2(t_5) \dot s_3 B^- 
$$
for $t_2,t_5\in \C^*$ and $m_4\in\C$ (see \eqref{e:DeoComp}), 
and they are not seen by $Z$.   

Since by Kostant \cite{Kos:QCoh} 
the full Peterson variety $Y$, 
and hence its open stratum $Y_B$, are irreducible,
this problem does not occur to the same extent 
in the full flag variety case. In that case the 
open embedding $Z^{\circ,crit}\to Y^\circ_B$ 
automatically has dense image, and for generic fixed 
value of $\tilde q=(\tilde q_1,\dotsc,\tilde q_n)$ the fiber of 
$Z^{\circ,crit}$ has the full number of 
of points (that is, $(n+1)!$).

\section{Total Positivity}\label{s:TotPos}
\subsection{}

The totally positive and nonnegative parts, 
$G_{>0}$ and $G_{\ge 0}$, of $G=SL_{n+1}$ are the semialgebraic subsets of 
$SL_{n+1}(\mathbb R)$ consisting of those matrices all of whose subdeterminants
are positive, respectively nonnegative. Equivalently, $g$ lies in 
$G_{\ge 0}$
if it acts by matrices with nonnegative real entries in all of the 
fundamental representations $\bigwedge^r \C^{n+1}$, with respect to 
the standard bases of these representations. Similarly 
$g$ belongs to $G_{>0}$ if it acts by matrices
with strictly positive entries. This
strong notion of positivity for $SL_{n+1}$, or
the general linear group, goes back to 
work of Schoenberg and Gantmacher and Krein
from the 1930's, see also \cite{Lus:IntroTotPos}. 

A useful characterization of $G_{\ge 0}$ is the following. Note that 
the simple root subgroups define semigroups  
$x_i(t),y_i(t)$ in $G_{\ge 0}$, where $t\in\R_{\ge 0}$. We also have
a semigroup given by the totally nonnegative part of the torus $T_{>0}$, 
the diagonal matrices with positive entries. 
By a theorem of Ann Whitney \cite{Whi:TotPosRed} these semigroups 
together generate $G_{\ge 0}$, and this description of $G_{\ge 0}$ was used by Lusztig  
 \cite{Lus:TotPos94} to extend the notion of total positivity to arbitrary reductive 
 algebraic groups. In fact, let $U^+_{\ge 0}$ and 
$U^-_{\ge 0}$ be the semigroups inside $U^+$ and $U^-$ generated
by the $\{x_i(t)\ |\ t\ge 0\}_{i\in I}$ and the $\{y_i(t)\ |\ t\ge 0\}_{i\in I}$, 
respectively. Then Lusztig noted that one has a `triangular decomposition'
\begin{equation*}
G_{\ge 0}=U^+_{\ge 0} T_{>0} U^-_{\ge 0},
\end{equation*} 
and also introduced a cell decomposition for $U^+_{\ge 0}$ 
-- and thereby for 
$U^-_{\ge 0}$ and $G_{\ge 0}$ -- which goes as follows. Let $w\in W$
and define 
\begin{equation*}
U^+(w):=U^+_{\ge 0}\cap B^-\dot w B^-.
\end{equation*}
If one chooses a reduced expression $\mathbf w=s_{i_1}\dotsc s_{i_m}$ 
for $w$, then $U^+(w)$ is shown to agree with the set
\begin{equation}\label{e:U(w)}
\{x_{i_1}(t_1)\dotsc x_{i_m}(t_m)\ |\ t_j\in\R_{>0}\},
\end{equation}
making it a semialgebraic cell of dimension $m$. The unique cell
of maximal dimension, $U^+(w_0)$, is also denoted
by $U^+_{>0}$.

Lusztig also defined a totally positive and a totally nonnegative part
for the flag variety $G/B^-$ (in our conventions), see  
\cite[Section~8]{Lus:TotPos94}.
These are given by
\begin{align*}
(G/B^-)_{>0}&:=
\{uB^- \ |\ u\in U^+_{>0}\},\\
(G/B^-)_{\ge 0}&:=\overline {(G/B^-)_{>0}},
\end{align*}
where the closure is taken inside the real flag variety $(G/B^-)(\R)$
with respect to its topology as a real manifold. By \cite{Rie:CelDec}
$(G/B^-)_{\ge 0}$ has a cell decomposition with cells 
$$
\mathcal R_{v,w;>0}:= (G/B^-)_{\ge 0}\cap \mathcal R_{v,w},
$$
as conjectured by Lusztig in \cite{Lus:IntroTotPos}.
An explicit description of these cells mimicking Lusztig's 
factorizations \eqref{e:U(w)} is the following 
\cite[Theorem~11.3]{MarRie:ansatz},
\begin{multline}
\label{e:Rpos}
\mathcal R_{ v, w;>0}=
\mathcal R_{\mathbf {v_+},\mathbf w}\cap \mathcal (G/B^-)_{\ge 0}\\
=\left\{g_1 g_2\dotsc g_mB_-\ \left | \ g_r=\begin{cases}
\dot s_{i_l}, 
& \text{if $r\in J^+_{(\mathbf {v_+},\mathbf w)}$,}\\
x_{i_r}(t_r),\ t_l\in\R_{>0}, &
\text{if $r\in J^\circ_{(\mathbf v,\mathbf w)}$}
\end{cases}\quad \right .\right\}.
\end{multline}
Here  
$\mathbf v_{+}$ 
is the positive subexpression $s_{i_{j_1}} s_{i_{j_2}}\dotsc
s_{i_{j_l}}$ for $v$ in  the reduced
expression $\mathbf w$ for $w$
from above, see also Section~\ref{s:Deodhar}.

We define the totally nonnegative parts of the Peterson variety and 
its strata by,
\begin{eqnarray*}
Y_{\ge 0}&:=&Y(\R)\cap (G/B^-)_{\ge 0},\\
Y_{P,\ge 0}&:=&Y_P(\R)\cap (G/B^-)_{\ge 0},\\
Y_{(P,P'),>0}&:=&Y_{(P,P')}(\R)\cap (G/B^-)_{\ge 0}=
Y(\R) \cap \mathcal R_{w_P,w_{P'};>0}. 
\end{eqnarray*}
The totally positive part of $Y$ is $Y_{>0}:=Y_{(B,G),>0}$.
\subsection{}
The GBCKS variety $Z_P$ 
also has a natural `positive part'. We set
\begin{align*}
Z_{P,\ge 0}&:=\{\rho\in Z_P \ |\ \rho_a\in\R_{\ge 0}\ \text{all}\
a\in\mathcal A\},\\
Z_{P,>0}&:=Z_{P,\ge 0}\cap Z_P^\circ,\\
Z_{(P,P'),>0}&:=Z_{(P,P')}\cap Z_{P,\ge 0}.
\end{align*}
Similarly, let $Z^{crit}_{P,\ge 0}:=Z_P^{crit}\cap Z_{\ge 0}$, and 
$Z^{crit}_{P,>0}:=Z_P^{crit}\cap Z_{P,> 0}$ and 
$Z^{crit}_{(P,P'),>0}:= Z_P^{crit}\cap Z_{(P,P'),>0}$.

\begin{prop}\label{p:pos} 
\begin{enumerate}
\item We have the following decomposition,
$$
Z^{crit}_{P,\ge 0}=\bigsqcup_{P'\supseteq P} Z^{crit}_{(P,P'),>0}.
$$
\item
The map $\phi^{crit}_P: Z^{crit}\to Y_P$
restricts to the positive strata giving homeomorphisms
$$\phi^{crit}_{(P,P'),>0}:
Z^{crit}_{(P,P'),>0}\overset\sim\To Y_{(P,P'),>0 }.
$$
\end{enumerate}
\end{prop}

\begin{proof}
To prove (1) it is sufficient to show that
\begin{equation*}Z^{crit}_{(P,P'),>0}=\{ \rho\in Z^{crit}_{P,\ge 0}\ | \ 
\text{$\tilde q_j(\rho)=0 \iff n_j\in I^{P'}$} \}.
\end{equation*} 
The inclusion ``$\subseteq$'' is clear. Let $I^{P'}=
\{n_{j_1},\dotsc,n_{j_t}\}\subset I^P$, with 
$1\le j_1<\dotsc < j_t\le k$. 
Suppose $\rho\in Z^{crit}_{P,\ge 0}$ with 
$\tilde q_{j_i}(\rho)=0$ for $i=1,\dotsc, t$, and
$\tilde q_l(\rho)\ne 0$ for all other $1\le l\le k$.

Let $v=(v_1,v_2)$ be a vertex in
$\mathcal V^{P'}_\bullet$. Then there is some $i$
such that $v_1\le n_{j_i}\le v_2$. Or in other 
words, the vertex $(n_{j_i},n_{j_i})$ in 
$\mathcal V_\bullet^{P'}$ 
lies above and to the right of $v$.
We need to show that any arrow $a$ for which either $h(a)=v$
or $t(a)=v$, satisfies $a(\rho)=0$. 

Recall the critical point condition at $v$,  
\begin{equation}\label{e:critQ}
\sum_{a', h(a')=v} a'(\rho)=\sum_{a'', t(a'')=v} a''(\rho).
\end{equation}
We suppose indirectly that one of these coordinates, either 
an $a'$ or $a''$, is nonzero 
on $\rho$. Let us call this coordinate $a_0$. 
Since $a(\rho)\ge 0$ for all $a\in \mathcal A$ 
it follows that both sides of the equation \eqref{e:critQ} must be nonzero.
So at least one of the coordinates on the opposite side of the equation 
to $a_0$ must also be positive on $\rho$. 
  
We can now define a sequence of coordinates, $a_{-m},a_{-m+1},\dotsc,
a_0,a_1,\dotsc, a_{m'}$, all of which should be nonzero 
on $\rho$, as follows. Start with $a_0$. If $a_i$ for $i\ge 0$ has 
been defined and has 
$h(a_i)\in\mathcal V^P_{\bullet}$, then
there is at least one arrow $a''$ with 
$t(a'')=h(a_i)$ and $a''(\rho)>0$. We set $a_{i+1}=a''$ (chosen 
arbitrarily if there are two such coordinates). 
The sequence ends when an arrow $a_{m'}$ has $h(a_{m'})
\in \mathcal V_\star^P$.     

On the other side, if $a_{-i}$ has been defined with $t(a_i)\in
\mathcal V_{\bullet}^P$, then from \eqref{e:critQ} it follows that
there is at least one $a'$ with $h(a')=t(a_{-i})$ and $a'(\rho)>0$.
So we set $a_{-i-1}=a'$. The sequence ends with $a_{-m}$ 
in the negative 
direction if $t(a_{-m})\in\mathcal V^P_\star$.  

Now $t(a_{-m})=\star_{l}$ and $h(a_{m'})=\star_{l'}$ where
$1\le l'\le l\le k+1$. By construction the 
vertex $\star_l$ is below and to the right of $v$, 
while the vertex $\star_{l'}$ is above and to the 
left of $v$. We have
$$
a_{-m}a_{-m+1}\dotsc a_{m'}=\tilde q_{l'} \tilde q_{l'+1}
\dotsc \tilde q_l.
$$
Since the vertex $(n_{j_i},n_{j_i})$ is above and to the right of $v$
it follows
that $l'\le j_i\le l$. Therefore the product
$\tilde q_{l'}\dotsc \tilde q_l$ vanishes on $\rho$
and we have the desired contradiction. 

Part (2) of the proposition
follows directly from the parameterization of 
the totally positive part of 
$\mathcal R_{w_P,w_{P'}}$ given in
\eqref{e:Rpos}, and also \cite{Lus:TotPos94} if $P=B$.
\end{proof}

\begin{thm}\label{t:theorem2}
Let $P'\supseteq P$ and 
$I^{P}\setminus I^{P'}=\{n_{k_1},\dotsc,n_{k_m}\}$. 
The restriction of the branched covering 
$q=(q_1,\dots,q_{k}):Y_P\to \C^{k}$
to the totally positive stratum $Y_{(P,P'),>0}$
gives rise to a homeomorphism
\begin{equation*}
(q_{k_1},\dotsc, q_{k_m}):Y_{(P,P'),>0}\To \R_{>0}^{|I^P|-|I^{P'}|}.
\end{equation*}
\end{thm}

\begin{proof}
By Proposition~\ref{p:pos} it is equivalent to show that 
$$
(\tilde q_{k_1},\dotsc, \tilde q_{k_m}):
Z^{crit}_{(P,P'),>0}\To \R_{>0}^{|I^P|-|I^{P'}|}
$$ 
is a homeomorphism. Assume first that 
$P'=G$ and let 
$\tilde Q\in\R_{> 0}^{k}$. Then the fiber 
$Z_{\tilde Q}$ lies inside $Z_P^\circ$
and we have to prove that 
$\mathcal F|_{Z_{\tilde Q}}$ has
a unique critical point 
in $Z_{\tilde Q}\cap Z_{P,>0}=:Z_{\tilde Q,>0}$. 

We begin by showing that a positive critical point (a minimum) exists. 
Since we are in $Z_{P,>0}$ we can write $\mathcal F$ in terms of 
the logarithmic 
vertex variables from Section~\ref{s:crit}. We have 
\begin{eqnarray*}
\R^\mathcal V&\To &Z_{P,>0}\\
(T_{v})_{v\in\mathcal V}&\mapsto &(e^{T_{h(a)}-T_{t(a)}})_{a\in \mathcal A}.
\end{eqnarray*}
Let us fix
$T_{\star_{j}}=T_{\star_{j}}(\tilde Q)=\sum_{i=j}^{k}\ln(\tilde Q_i)$
and  $T_{\star_{k+1}}=T_{\star_{k+1}}(\tilde Q)=0$.
Then the above map restricts to a diffeomorphism
\begin{equation}\label{e:diffeo}
\R^{\mathcal V_\bullet}\overset\sim\To Z_{\tilde Q,>0}.
\end{equation}
We now define $\mathcal F_{\tilde Q,>0}$ to be
the restriction of $\mathcal F$ to $Z_{\tilde Q,>0}$ and 
identify $Z_{\tilde Q,>0}$ with $\R^{\mathcal V_{\bullet}}$
by \eqref{e:diffeo}. 
So 
\begin{eqnarray*}
\mathcal F_{\tilde Q,>0}:\R^{\mathcal V_{\bullet}}&\to&\R_{>0}\\
\qquad \quad (T_{v})_{v\in\mathcal V_{\bullet}}&\mapsto&
\sum_{a\in \mathcal A}e^{T_{h(a)}-T_{t(a)}}.
\end{eqnarray*}
Define nested subsets 
$C_1\subset C_2\subset C_3\dotsc \subset \R^{\mathcal V_{\bullet}} $ by
$$
C_m:=\{(T_v)\in \R^{\mathcal V_\bullet}\ 
|\ \text{ $ T_{\star_1}(\tilde Q)- m |\mathcal A|\le T_v\le
T_{\star_{k+1}}(\tilde Q)+ m |\mathcal A|$, for all $v\in
\mathcal V_\bullet$ }\}.
$$    
The set $C_m$  is nonempty 
for sufficiently large $m$, and clearly compact. 

Let $(\tilde T_v)$ be an element of $Z_{\tilde Q,>0}$ 
not in $C_m$. Then there are two possibilities
\begin{enumerate}
\item
Suppose first that there is a vertex $v\in\mathcal V_\bullet$ 
such that 
$$
T_{\star_1}(\tilde Q)-\tilde T_v>m |\mathcal A|.
$$
We can find a sequence of vertices  $v_0=v,v_1,\dotsc, v_s=\star_1$ and
a sequence of arrows $a_1,\dotsc, a_s$ such that $t(a_i)=v_{i-1}$ and
$h(a_i)=v_i$. So we have 
$$
T_{\star_1}(\tilde Q)-\tilde T_v=\sum_{i=1}^s (\tilde T_{h({a_i})}-\tilde 
T_{t({a_i})})> m |\mathcal A|.
$$
Since there are fewer 
than $|\mathcal A|$ summands, one of the summands must satisfy 
$\tilde T_{h({a_i})}-\tilde T_{t({a_i})}>m$. Therefore 
we have 
$$\mathcal F((\tilde T_v))
>e^{\tilde T_{h(a_i)}-\tilde T_{t({a_i})}}>e^m.
$$
\item
Otherwise we have $\tilde T_{v}-T_{\star_{k+1}}(\tilde Q)>m|\mathcal A|$
for some vertex $v\in\mathcal V_{\bullet}$. 
In this case the analogous argument to above, but using a path 
from $\star_{k+1}$ to $v$, implies that 
$\mathcal F((\tilde T_v))>e^m$. 
\end{enumerate}

Now for every $m$ with $C_m$ nonempty,
$\mathcal F_{\tilde Q,>0}|_{C_m}$ attains a minimum $c_m$. 
And the sequence of minima, $c_m\ge c_{m+1}\ge \cdots$, 
stabilizes to give a global minimum $c$ for $\mathcal F_{\tilde Q,>0}$,
since eventually $c_{m'}< e^{m'}$.  

Let us now show 
that there are no other critical points 
in $Z_{\tilde Q,> 0}$. For this it suffices to prove that the 
Hessian of $\mathcal F_{\tilde Q,> 0}$ is everywhere positive definite,
which follows by direct calculation~:

\begin{multline*}
\left
(\sum_{v\in\mathcal V_\bullet} m_v\frac{\partial}{\partial T_v}
\right)^2 \mathcal F\\
=\left(\sum_{v\in\mathcal V_\bullet} m_v\frac{\partial}{\partial T_v} 
\right)
\sum_{v'\in \mathcal V_\bullet} m_{v'}\left(\sum_{a,\, h(a)=v'} 
e^{T_{v'}-T_{t(a)}} -
\sum_{a,\, t(a)=v'} e^{T_{h(a)}-T_{v'}}\right)=\\
\sum_{v\in\mathcal V_{\bullet}} 
m_v^2\left(\sum_{a,\, h(a)=v} e^{T_{v}- T_{t(a)}} +
\sum_{a,\, t(a)=v} e^{T_{h(a)}-T_v}\right) -
2\sum_{a\in\mathcal A}m_{h(a)}m_{t(a)} e^{T_{h(a)}-T_{t(a)}}=\\
\sum_{a\in\mathcal A}(m_{h(a)}-m_{t(a)})^2 e^{T_{h(a)}-T_{t(a)}}.
\end{multline*}

Viewing $\mathcal F_{\tilde Q,>0}$ as a family of functions
on $\R^{\mathcal V_{\bullet}}$, we have shown that 
for each $\tilde Q$ there is a unique 
$\rho_{\tilde Q}\in \R^{\mathcal V_{\bullet}}$ such that
$$
\mathcal F_{\tilde Q,>0}(\rho_{\tilde Q})=c_{\tilde Q}
$$ 
is a minimum for $\mathcal F_{\tilde Q,>0}$. 
Since $\mathcal F_{\tilde Q,>0}$ depends
continuously on $\tilde Q$ it 
follows that so does the point $\rho_{\tilde Q}$.
Therefore
$(\tilde q_1,\dotsc, \tilde q_k):Z^{crit}_{(P,G),>0}\to \R_{>0}^k$ 
is a homeomorphism.

In the case of $Z_{(P,P'),>0}$  for general $P'\supseteq P$
the graph $(\mathcal V,\mathcal A)$ can be replaced by 
a disjoint union of subgraphs, in each of which all edges correspond to
strictly 
nonzero coordinates on $Z_{(P,P')}$. In this case the same arguments 
as above,
now applied to each one of the
subgraphs, prove the theorem.
\end{proof}

\subsection{The cell decomposition of $Y_{\ge 0}$.}
Theorem~\ref{t:theorem2} gives a cell decomposition of 
$Y_{P,\ge 0}$ for every $P$. 
Therefore in total we have a cell decomposition for all of $Y_{\ge 0}$
with  cells indexed by pairs of parabolics  $(P,P')$ satisfying
$B\subseteq P\subseteq P'$.

Recall that $I=\{1,\dotsc, n\}$. Consider the set $\mathcal J$ 
of pairs $(J,K)$ of subsets $J,K\subseteq I$ with $J\subseteq K$.
This set is a poset under the partial ordering
 $$
(J,K)\le (J',K') :\iff  J'\subseteq J\subseteq K\subseteq K'. 
 $$
Moreover $\mathcal J$ can be identified with the face poset
of the $n$-dimensional cube $[0,1]^I$ by
\begin{equation*}
F_{(J,K)}=\left\{(x_i)\in [0,1]^I\ \left |\ \begin{array}{l}
\text{$x_i=0$ if $i\in J$,}\\
\text{$0<x_i<1$ if $i\in K\setminus J$,}\\
\text{$x_i=1$ if $i\notin K$}
\end{array}\right.\right \}
\end{equation*}
That is, $F_{(J,K)}$ is an open face of
the hypercube $[0,1]^I$ 
of dimension $|K\setminus J|$. 

\begin{conj}\label{c:cube}
There is a homeomorphism
\begin{equation*}
Y_{\ge 0}\to [0,1]^I
\end{equation*}
such that $Y_{(P,P'),>0}$ is mapped to $F_{(I^{P'},I^{P})}$.  
\end{conj}

Recall that in \cite[Corollary~7.4]{Rie:QCohPFl} we constructed a 
homeomorphism 
\begin{equation}\label{e:TotPosToepl}
Y_{\ge 0}\cap B^-\dot w_0 B^-/B^-=
\bigsqcup_P Y_{(P,G),>0}\overset\sim\To\R^{n}_{\ge 0}
\end{equation}
of cell decomposed spaces. 
Therefore $Y_{\ge 0}$ indeed resembles a cube
in a neighborhood of the fixed point $Y_{(G,G)}=\{\dot w_0 B^-\}$.

We end by proving two corollaries of Theorem~\ref{t:theorem2} 
which give further evidence for Conjecture~\ref{c:cube}. 
Firstly we see that the 
parameterizations from Theorem~\ref{t:theorem2} can be 
combined similarly to \eqref{e:TotPosToepl} above. In particular  
$Y_{\ge 0}$ also looks like a cube
in a neighborhood of $Y_{(B,B)}=\{B^-\}$. 
\begin{cor}\label{c:finite}
The map $q=(q_1,\dotsc, q_k):Y_{P,\ge 0}\To \R_{\ge 0}^{k}$ 
is a homeomorphism. 
\end{cor}
\begin{proof}
Since $qH^*(G/P)=H^*(G/P)\otimes\C[q_1,\dotsc, q_k] $ it follows that
the map 
$$(q_1,\dotsc, q_k): Y_P\to \C^k$$ 
is finite. Therefore 
the restriction to the closed subset $Y_{P,\ge 0}$ of $Y_P$ is
proper and in particular closed.
By Theorem~\ref{t:theorem2}  the map
$(q_1,\dotsc, q_k):Y_{P,\ge 0}\to \R^k_{\ge 0}$ is a 
bijection. Since it is also
continuous and closed it is a homeomorphism.
\end{proof}

Finally, we can use this result to show 
that $Y_{\ge 0}$ is contractible.
\begin{cor}
The totally nonnegative part of the Peterson
variety is contractible. 
\end{cor}

\begin{proof}
We claim first that any element of
$U^-_{>0}$ translates the
totally nonnegative part of $G/B^-$
into the totally nonnegative
part of the big cell $B^+ B^-/B^-$. 
This can be proved one (opposite) Bruhat 
cell at a time. Let 
us consider the totally nonnegative part of $B^+\dot w B^-/B^-$
and act on it by some totally positive
$u\in U^-$.
Indeed, in this case we can factorize $u\in U^-(w_0)$ into
$u=u_1 u_2$ with $u_1\in U^-(w_0w)$ and
$u_2\in U^-(w\inv)$. Then using $u_2\in 
B^+ \dot w\inv B^+$ and the properties of Bruhat decomposition
we see that
\begin{multline*}
u_2 B^+\dot w B^-/B^-\quad \subseteq\quad
(B^+ \dot w\inv B^+)\dot w\ (\dot w_0 B^+ \dot w_0\inv)\ B^-/B^-\\
=(B^+ \dot w\inv B^+ \dot w \dot w_0 B^+) \dot w_0\inv B^-/B^-=
B^+ \dot w_0 B^+\dot w_0\inv B^-/B^-= B^+ B^-/B^-.
\end{multline*}
Since by \cite[Proposition~8.13]{Lus:TotPos94}  $G_{\ge 0}$
preserves the totally nonnegative part of the flag variety,
it follows that
$$
u (B^+ \dot w B^-/B^-)_{\ge 0}= u_1 u_2 (B^+\dot w B^-/B^-)_{\ge 0}
=u_1 (B^+ B^-/B^-)_{\ge 0} = u_1 U^+_{\ge 0} B^-/B^-.
$$
Now it follows from \cite[Lemma~2.3]{Lus:TotPos94} that 
$u_1U^+_{\ge 0}\subset U^+_{\ge 0}T_{>0}U^-_{\ge 0}$.
Therefore we have
$u (B^+\dot w B^-/B^-)_{\ge 0}\subset (B^+B^-/B^-)_{\ge 0}$ as 
required.

Let $y(t):=\exp(t f)$ for $t\in \R_{\ge 0}$. Then
$y(t)\in U^-_{>0}$, 
for $t>0$, by \cite[Proposition~5.9]{Lus:TotPos94}. 
Therefore we have
$y(t) (G/B^-)_{\ge 0}\subseteq
(B^+B^-/B^-)_{\ge 0}$ for all $t>0$.  
Since the action of $y(t)$ on the 
flag variety preserves the 
Peterson variety this implies
\begin{eqnarray*}
  y(t)\cdot Y_{\ge 0} \subset Y_{B,\ge 0},& \text{ for all $t>0$.}
\end{eqnarray*}
By Corollary~\ref{c:finite}, $Y_{B,\ge 0}$ is contractible. 
Let $F'_s$ be a deformation retraction 
$F'_s:Y_{B,\ge 0}\to Y_{B,\ge 0}$ such that $F'_0=id$ and 
$F'_1(Y_{B,\ge 0})=\{B^-\}$, and $s\in [0,1]$. 
Then 
\begin{equation*}
F_s(gB^-):=\begin{cases}
gB^- &\text{ if $s=0$},\\
F'_s(y(s)gB^-) &\text{ if $0<s\le 1$,} 
\end{cases}
\end{equation*}
defines a deformation retraction $F_s:Y_{\ge 0}\to Y_{\ge 0}$ with the 
same properties, 
implying that $Y_{\ge 0}$ is also contractible. 
\end{proof}

\bibliographystyle{amsplain}

\def\cprime{$'$}
\providecommand{\bysame}{\leavevmode\hbox to3em{\hrulefill}\thinspace}
\providecommand{\MR}{\relax\ifhmode\unskip\space\fi MR }
\providecommand{\MRhref}[2]{%
  \href{http://www.ams.org/mathscinet-getitem?mr=#1}{#2}
}
\providecommand{\href}[2]{#2}

\end{document}